\documentclass[11pt,letterpaper,colorlinks=true,linkcolor=blue,citecolor=blue,urlcolor=blue]{scrartcl}
\usepackage[T1]{fontenc}
\usepackage[utf8]{inputenc}
\usepackage{lmodern}
\usepackage[english]{babel}
\usepackage[english]{isodate}
\usepackage{amsthm}
\usepackage{amsmath}
\usepackage{amssymb}
\usepackage{amsfonts}
\usepackage{mathtools}
\usepackage{mathrsfs}
\usepackage{dsfont}
\usepackage{upgreek}
\usepackage{xpatch}
\usepackage{bm}
\usepackage{bbm}
\usepackage{enumitem}
\usepackage[toc,page]{appendix}
\usepackage{footnotebackref}
\usepackage[totalwidth=480pt, totalheight=680pt]{geometry}
\usepackage{hyperref}

\theoremstyle{remark}
\newtheorem{assumptions}{Assumptions}[section]
\newtheorem{remark}[assumptions]{Remark}

\theoremstyle{plain}
\newtheorem{theorem}[assumptions]{Theorem}
\newtheorem{proposition}[assumptions]{Proposition}
\newtheorem{lemma}[assumptions]{Lemma}
\newtheorem{corollary}[assumptions]{Corollary}

\addtokomafont{disposition}{\normalfont\scshape}

\makeatletter
\xpatchcmd{\@sec@pppage}{
\bfseries}{
\normalfont\scshape\Large}{}{}
\makeatother

\setkomafont{section}{\large}
\setkomafont{subsection}{\normalsize}

\numberwithin{equation}{section}

\begin{document}


\title{
\LARGE A variational characterization of Langevin--Smoluchowski diffusions\thanks{
We are indebted to Walter Schachermayer for joint work that this paper complements, and we are grateful to Promit Ghosal for his interest in this work and for encouraging us to consider its iterative aspects. We would also like to thank Michele Pavon for pointing out relevant literature to us.

\smallskip

\noindent I. Karatzas acknowledges support from the National Science Foundation (NSF) under grant NSF-DMS-20-04997. \newline
B. Tschiderer acknowledges support by the Austrian Science Fund (FWF) under grant P28661 and by the Vienna Science and Technology Fund (WWTF) through project MA16-021.
}
\vspace{0.6cm}
}

\subtitle{
\normalsize To the memory of Mark H.A.\ Davis, 1945 -- 2020
\vspace{0.2cm}
}

\author{
\large Ioannis Karatzas\thanks{
Department of Mathematics, Columbia University, 2990 Broadway, New York, NY 10027, USA 
\newline (email: \href{mailto: ik1@columbia.edu}{ik1@columbia.edu}); 
\newline and INTECH Investment Management, One Palmer Square, Suite 441, Princeton, NJ 08542, USA 
\newline (email: \href{mailto: ikaratzas@intechjanus.com}{ikaratzas@intechjanus.com}).
} 
\and \large Bertram Tschiderer\thanks{
Faculty of Mathematics, University of Vienna, Oskar-Morgenstern-Platz 1, 1090 Vienna, Austria \newline (email: \href{mailto: bertram.tschiderer@univie.ac.at}{bertram.tschiderer@univie.ac.at}).
}
}


\date{\normalsize 7th November 2020}


\maketitle


\begin{abstract} \small \noindent \textsc{Abstract.} We show that Langevin--Smoluchowski measure on path space is invariant under time-reversal, followed by stochastic control of the drift with a novel entropic-type criterion. Repeated application of these forward-backward steps leads to a sequence of stochastic control problems, whose initial/terminal distributions converge to the Gibbs probability measure of the diffusion, and whose values decrease to zero along the relative entropy of the Langevin--Smoluchowski flow with respect to this Gibbs measure.

\bigskip

\small \noindent \href{https://mathscinet.ams.org/mathscinet/msc/msc2010.html}{\textit{MSC 2010 subject classifications:}} Primary 93E20, 94A17; secondary 35Q84, 60G44, 60J60 

\bigskip

\end{abstract}


\section{Introduction} \label{sec.one}


Diffusions of Langevin--Smoluchowski type have some important properties. They possess invariant (Gibbs) probability measures described very directly in terms of their potentials and towards which, under appropriate conditions, their time-marginals converge as time increases to infinity and in a manner that conforms to the second law of thermodynamics: the relative entropy of the current distribution, with respect to the invariant one, decreases to zero. The seminal paper \cite{JKO98} revealed another remarkable, local aspect of this decrease towards equilibrium: the family of time-marginals is, at (almost) every point in time, a curve of steepest descent among all probability density functions with finite second moment, when distances are measured according to the Wasserstein metric in configuration space.

\smallskip
 
We establish in this paper yet another variational property of such diffusions, this time a global one: their law is invariant under the combined effects of time-reversal, and of stochastic control of the drift under a novel, entropic-type criterion. Here, one minimizes over admissible controls the relative entropy of the ``terminal'' state with respect to the invariant measure, plus an additional term thought of as ``entropic cost of time-reversal'': the difference in relative entropy with respect to the Langevin--Smoluchowski measure on path space, computed based on the ``terminal'' state as opposed to on the entire path. Quite similar, but different, cost criteria have been considered in \cite{DPP89,DPP90,FS85,Pav89,PW91}.

\smallskip

The setting under consideration bears similarities to the celebrated Schr\"odinger bridge problem, but also considerable differences. Both problems are posed on a fixed time horizon of finite length, and both involve the relative entropy with respect to the invariant measure. But here this entropy is modified by the addition of the above-mentioned entropic cost of time-reversal, and there is no fixed, target distribution on the terminal state. Yet the trajectory that emerges as the solution of the stochastic control problem has time-marginals that replicate exactly those of the original Langevin--Smoluchowski flow, whence the ``invariance'' property mentioned in the abstract. 

\smallskip

We refer to \cite{CGP16a,CGP16b,CGP20,Leo14a} for overviews on the classical Schr\"odinger bridge problem, to \cite{TT13} for the related semimartingale transport problem, and to the recent paper \cite{BCGL20} for a detailed study of the mean-field Schr\"odinger problem. A related controllability problem for a Fokker--Planck equation and its connection to Schr\"odinger systems and stochastic optimal control, is considered in \cite{Bla92}. More information about the Schr\"odinger equation, diffusion theory, and time reversal, can be found in the book \cite{Nag93}. 


\subsection{Preview}


In \hyperref[the.setting]{Section \ref*{the.setting}} we introduce the Langevin--Smoluchowski measure $\mathds{P}$ on path space, under which the canonical process $(X(t))_{t \geqslant 0}$ has dynamics \hyperref[ls.dp.sde]{(\ref*{ls.dp.sde})} with initial distribution $P(0)$. Then, in \hyperref[sec.three]{Section \ref*{sec.three}}, this process is studied under time-reversal. That is, we fix a terminal time $T \in (0,\infty)$ and consider the time-reversed process $\overline{X}(s) = X(T-s)$, $0 \leqslant s \leqslant T$. Classical theory of time reversal shows that $\overline{X}$ is again a diffusion, and gives an explicit description of its dynamics.

\smallskip

\hyperref[sec.ascp.f]{Section \ref*{sec.ascp.f}} develops our main result, \hyperref[m.theo.one]{Theorem \ref*{m.theo.one}}. An equivalent change of probability measure $\mathds{P}^{\boldsymbol{\gamma}} \sim \mathds{P}$ adds to the drift of $\overline{X}$ a measurable, adapted process $\boldsymbol{\gamma}(T-s)$, $0 \leqslant s \leqslant T$. In broad brushes, this allows us to define, in terms of relative entropies, the quantities
\begin{equation} \label{inf.inf.hd}
H^{\boldsymbol{\gamma}} \vcentcolon = H\big( \mathds{P}^{\boldsymbol{\gamma}} \, \vert \, \mathds{Q}\big) \Big\vert_{ \boldsymbol{\sigma}(\overline{X}(T))}
\, , \qquad
D^{\boldsymbol{\gamma}} \vcentcolon = H\big( \mathds{P}^{\boldsymbol{\gamma}} \, \vert \, \mathds{P}\big) \Big\vert_{ \boldsymbol{\sigma}(\overline{X})} 
- H\big( \mathds{P}^{\boldsymbol{\gamma}} \, \vert \, \mathds{P}\big) \Big\vert_{\boldsymbol{\sigma}(\overline{X}(T))}.
\end{equation}
Here, $\mathds{Q}$ is the probability measure on path space, inherited from the Langevin--Smoluchowski dynamics \hyperref[ls.dp.sde]{(\ref*{ls.dp.sde})} with initial distribution given by the invariant Gibbs probability measure $\mathrm{Q}$. \hyperref[m.theo.one]{Theorem \ref*{m.theo.one}} establishes then the variational characterization
\begin{equation} \label{inf.inf}
\inf_{\boldsymbol{\gamma}} \big( H^{\boldsymbol{\gamma}} + D^{\boldsymbol{\gamma}} \big) = H\big( P(T) \, \vert \, \mathrm{Q}\big),
\end{equation}
where $P(T)$ denotes the distribution of the random variable $X(T)$ under $\mathds{P}$. The process $\boldsymbol{\gamma}_{\ast}$ that realizes the infimum in \hyperref[inf.inf]{(\ref*{inf.inf})} gives rise to a probability measure $\mathds{P}^{\boldsymbol{\gamma}_{\ast}}$, under which the time-reversed diffusion $\overline{X}$ is of Langevin--Smoluchowski type in its own right, but now with initial distribution $P(T)$. In other words, with the constraint of minimizing the sum of the entropic quantities $H^{\boldsymbol{\gamma}}$ and $D^{\boldsymbol{\gamma}}$ of \hyperref[inf.inf.hd]{(\ref*{inf.inf.hd})}, Langevin--Smoluchowski measure on path space is invariant under time-reversal.

\smallskip

\hyperref[sec.ltm]{Sections \ref*{sec.ltm}} -- \hyperref[conclusion.seven]{\ref*{conclusion.seven}} develop ramifications of the main result, including the following consistency property: starting with the time-reversal $\overline{X}$ of the Langevin--Smoluchowski diffusion, the solution of a related optimization problem, whose value is now $H(P(2T) \, \vert \, \mathrm{Q})$, leads to the original forward Langevin--Smoluchowski dynamics, but now with initial distribution $P(2T)$. Iterating these procedures we obtain an alternating sequence of forward-backward Langevin--Smoluchowski dynamics with initial distributions $(P(kT))_{k \in \mathds{N}_{0}}$ converging to $\mathrm{Q}$ in total variation, along which the values of the corresponding optimization problems as in \hyperref[inf.inf]{(\ref*{inf.inf})} are given by $(H(P(kT) \, \vert \, \mathrm{Q}))_{k \in \mathds{N}}$ and decrease to zero.


\section{The setting} \label{the.setting}


Let us consider a Langevin--Smoluchowski diffusion process $(X(t))_{t \geqslant 0}$ of the form
\begin{equation} \label{ls.dp.sde}
\mathrm{d}X(t) = - \nabla \Psi\big(X(t)\big) \, \mathrm{d}t + \mathrm{d}W(t),
\end{equation}
with values in $\mathds{R}^{n}$. Here $(W(t))_{t \geqslant 0}$ is standard $n$-dimensional Brownian motion, and the ``potential'' $\Psi \colon \mathds{R}^{n} \rightarrow [0,\infty)$ is a $C^{\infty}$-function growing, along with its derivatives of all orders, at most exponentially as $\vert x \vert \rightarrow \infty$; we stress that no convexity assumptions are imposed on this potential. We posit also an ``initial condition'' $X(0) = \Xi$, a random variable independent of the driving Brownian motion and with given distribution $P(0)$. For concreteness, we shall assume that this initial distribution has a continuous probability density function $p_{0}(\, \cdot \,)$.

\smallskip

Under these conditions, the Langevin--Smoluchowski equation \hyperref[ls.dp.sde]{(\ref*{ls.dp.sde})} admits a pathwise unique, strong solution, up until an ``explosion time'' $\mathfrak{e}$; such explosion never happens, i.e., $\mathds{P}(\mathfrak{e}= \infty) = 1$, if in addition the second moment condition \hyperref[int.req.a]{(\ref*{int.req.a})} and the coercivity condition \hyperref[docc]{(\ref*{docc})} below hold. The condition \hyperref[docc]{(\ref*{docc})} propagates the finiteness of the second moment to the entire collection of time-marginal distributions $P(t) = \mathrm{Law}(X(t))$, $t \geqslant 0$, which are then determined uniquely. In fact, adapting the arguments in \cite{Rog85} to the present situation, we check that each time-marginal distribution $P(t)$ has probability density $p(t, \, \cdot \,)$ such that the resulting function $(t,x) \mapsto p(t,x)$ is continuous and strictly positive on $(0,\infty) \times \mathds{R}^{n}$; differentiable with respect to the temporal variable $t$ for each $x \in \mathds{R}^{n}$; smooth in the spatial variable $x$ for each $t > 0$; and such that the logarithmic derivative $(t,x) \mapsto \nabla \log p(t,x)$ is continuous on $(0,\infty) \times \mathds{R}^{n}$. These arguments also lead to the \textit{Fokker--Planck} \cite{Fri75,Gar09,Ris96,Sch80}, or \textit{forward Kolmogorov} \cite{Kol31}, equation 
\begin{equation} \label{fp.fk.pde}
\partial p(t,x) = \tfrac{1}{2} \Delta p(t,x) + \operatorname{div}\big(\nabla \Psi(x) \, p(t,x) \big), \qquad (t,x) \in (0,\infty) \times \mathds{R}^{n}
\end{equation}
with initial condition
\begin{equation}
p(0,x) = p_{0}(x), \qquad x \in \mathds{R}^{n}.
\end{equation}
Here and throughout this paper, $\partial$ denotes differentiation with respect to the temporal argument; whereas $\nabla$, $\Delta$ and $\operatorname{div}$ stand, respectively, for gradient, Laplacian and divergence with respect to the spatial argument.


\subsection{Invariant measure, likelihood ratio, and relative entropy} \label{subsec.2.1}


We introduce now the function 
\begin{equation} \label{dfq}
q(x) \vcentcolon = \mathrm{e}^{-2 \Psi(x)}, \qquad x \in \mathds{R}^{n}
\end{equation}
and note that it satisfies the stationary version
\begin{equation} \label{fp.fk.pde.st}
\tfrac{1}{2} \Delta q(x) + \operatorname{div}\big(\nabla \Psi(x) \, q(x) \big) = 0, \qquad x \in \mathds{R}^{n}
\end{equation}
of the forward Kolmogorov equation \hyperref[fp.fk.pde]{(\ref*{fp.fk.pde})}. We introduce also the $\sigma$-finite measure $\mathrm{Q}$ on the Borel subsets $\mathscr{B}(\mathds{R}^{n})$ of $\mathds{R}^{n}$, which has density $q$ as in \hyperref[dfq]{(\ref*{dfq})} with respect to $n$-dimensional Lebesgue measure. This measure $\mathrm{Q}$ is invariant for the diffusion of \hyperref[ls.dp.sde]{(\ref*{ls.dp.sde})}; see Exercise 5.6.18 in \cite{KS88}. When finite, $\mathrm{Q}$ can be normalized to an invariant probability measure for \hyperref[ls.dp.sde]{(\ref*{ls.dp.sde})}, to which the time-marginals $P(t)$ ``converge'' as $t \rightarrow \infty$; more about this convergence can be found in \hyperref[conclusion.seven]{Section \ref*{conclusion.seven}}. We shall always assume tacitly that such a normalization has taken place, when $\mathrm{Q}$ is finite, i.e., when
\begin{equation} \label{prob.gib}
\mathrm{Q}(\mathds{R}^{n}) = \int_{\mathds{R}^{n}} q(x) \, \mathrm{d}x = \int_{\mathds{R}^{n}} \mathrm{e}^{-2\Psi(x)} \, \mathrm{d}x < \infty.
\end{equation}
 
\smallskip

One way to think of the above-mentioned convergence, is by considering the likelihood ratio 
\begin{equation} \label{l.r.f}
\ell(t,x) \vcentcolon = \frac{p(t,x)}{q(x)}, \qquad (t,x) \in (0,\infty) \times \mathds{R}^{n}.
\end{equation}
It follows from \hyperref[fp.fk.pde]{(\ref*{fp.fk.pde})}, \hyperref[fp.fk.pde.st]{(\ref*{fp.fk.pde.st})} that this function satisfies the \textit{backward} Kolmogorov equation
\begin{equation} \label{b.k.pde}
\partial \ell(t,x) = \tfrac{1}{2} \Delta \ell(t,x) - \big\langle \nabla \ell(t,x) \, , \nabla \Psi(x) \big\rangle, \qquad (t,x) \in (0,\infty) \times \mathds{R}^{n}. 
\end{equation}

\smallskip

In terms of the likelihood ratio function \hyperref[l.r.f]{(\ref*{l.r.f})}, let us consider now the relative entropy
\begin{equation} \label{relative.entropy}
H\big( P(t) \, \vert \, \mathrm{Q} \big) \vcentcolon 
= \mathds{E}_{\mathds{P}}\big[ \log \ell\big(t,X(t)\big) \big]
= \int_{\mathds{R}^{n}} \log \bigg(\frac{p(t,x)}{q(x)}\bigg) \, p(t,x) \, \mathrm{d}x , \qquad t \geqslant 0
\end{equation}
of the probability distribution $P(t)$ with respect to the invariant measure $\mathrm{Q}$. The expectation in \hyperref[relative.entropy]{(\ref*{relative.entropy})} is well-defined in $[0,\infty]$, if $\mathrm{Q}$ is a probability measure. As we are not imposing this as a blanket assumption, we shall rely on Appendix C of \cite{KST20b}, where it is shown that the relative entropy $H( P(t) \, \vert \, \mathrm{Q})$ is well-defined and takes values in $(-\infty,\infty]$, whenever $P(t)$ belongs to the space $\mathscr{P}_{2}(\mathds{R}^{n})$ of probability measures with finite second moment (see also \cite{Leo14b} for a more general discussion). This, in turn, is the case whenever $P(0)$ has finite second moment, and the coercivity condition
\begin{equation} \label{docc}
\forall \, x \in \mathds{R}^{n}, \vert x \vert \geqslant R \colon \qquad \big\langle x \, , \nabla \Psi(x) \big\rangle \geqslant - c \, \vert x \vert^{2}
\end{equation}
is satisfied by the potential $\Psi$ in \hyperref[ls.dp.sde]{(\ref*{ls.dp.sde})}, for some real constants $c \geqslant 0$ and $R \geqslant 0$; see the first problem on page 125 of \cite{Fri75}, or Appendix B in \cite{KST20b}. The prototypical such potential is $\Psi(x) = \frac{1}{2} \vert x \vert^{2}$, leading to Ornstein--Uhlenbeck dynamics in \hyperref[ls.dp.sde]{(\ref*{ls.dp.sde})}; but $\Psi \equiv 0$ and the ``double well'' $\Psi(x) = (x^{2} - a^{2})^{2}$ for $a > 0$, are also allowed. In particular, the coercivity condition \hyperref[docc]{(\ref*{docc})} does not imply that the potential $\Psi$ is convex.

\smallskip

\textit{We shall impose throughout \textnormal{\hyperref[the.setting]{Sections \ref*{the.setting}} -- \hyperref[yet.st.cont.]{\ref*{yet.st.cont.}}} the coercivity condition \textnormal{\hyperref[docc]{(\ref*{docc})}}, as well as the finite second moment condition}
\begin{equation} \label{int.req.a}
\int_{\mathds{R}^{n}} \vert x \vert^{2} \, p_{0}(x) \, \mathrm{d}x < \infty.
\end{equation}
This amounts to $P(0) \in \mathscr{P}_{2}(\mathds{R}^{n})$, as has been already alluded to. In \hyperref[conclusion.seven]{Section \ref*{conclusion.seven}} we will see that these two conditions \hyperref[docc]{(\ref*{docc})} and \hyperref[int.req.a]{(\ref*{int.req.a})} are not needed when $\mathrm{Q}$ is a probability measure. 

\smallskip

However, \textit{we shall impose throughout the entire paper the crucial assumption that the initial relative entropy is finite}, i.e., 
\begin{equation} \label{int.req.b}
H\big( P(0) \, \vert \, \mathrm{Q}\big) = \int_{\mathds{R}^{n}} \log \bigg(\frac{p_{0}(x)}{q(x)}\bigg) \, p_{0}(x) \, \mathrm{d}x < \infty;
\end{equation}
under either the conditions ``\hyperref[docc]{(\ref*{docc})} $+$ \hyperref[int.req.a]{(\ref*{int.req.a})}'', or the condition \hyperref[prob.gib]{(\ref*{prob.gib})}, the decrease of the relative entropy\footnote{A classical aspect of thermodynamics; for a proof of this fact under the conditions ``\hyperref[docc]{(\ref*{docc})} $+$ \hyperref[int.req.a]{(\ref*{int.req.a})}'' and without assuming finiteness of $\mathrm{Q}$, see Theorem 3.1 in \cite{KST20a}; when $\mathrm{Q}$ is a probability measure, we refer to \hyperref[app.a.a]{Appendix \ref*{app.a.a}}.} function $[0,\infty) \ni t \mapsto H( P(t) \, \vert \, \mathrm{Q}) \in (-\infty,\infty]$ implies then that the quantity $H(P(t) \, \vert \, \mathrm{Q})$ in \hyperref[relative.entropy]{(\ref*{relative.entropy})} is finite for all $t \geqslant 0$ whenever \hyperref[int.req.b]{(\ref*{int.req.b})} holds.

\smallskip

In fact, under the conditions \hyperref[docc]{(\ref*{docc})} -- \hyperref[int.req.b]{(\ref*{int.req.b})}, the rate of decrease for the relative entropy, measured with respect to distances traveled in $\mathscr{P}_{2}(\mathds{R}^{n})$ in terms of the quadratic Wasserstein metric
\[
W_{2}(\mu,\nu) = \Big( \, \inf_{ Y \sim \mu, Z \sim \nu} \mathds{E} \vert Y - Z \vert^{2} \,  \Big)^{1/2} \, , \qquad \mu, \nu \in \mathscr{P}_{2}(\mathds{R}^{n})
\]
(cf.\ \cite{AG13,AGS08,Vil03}) is, at Lebesgue-almost all times $t_{0} \in [0, \infty)$, the \textit{steepest possible} along the Langevin--Smoluchowski curve $(P(t))_{t \geqslant 0}$ of probability measures. Here, we are comparing the curve $(P(t))_{t \geqslant 0}$ against all such curves $(P^{\beta}(t))_{t \geqslant t_{0}}$ of probability measures generated as in \hyperref[ls.dp.sde]{(\ref*{ls.dp.sde})} --- but with an additional drift $\nabla B$ for suitable (smooth and compactly supported) perturbations $B$ of the potential $\Psi$ in \hyperref[ls.dp.sde]{(\ref*{ls.dp.sde})}. This local optimality property of Langevin--Smoluchowski diffusions is due to \cite{JKO98}; it was established by \cite{KST20a} in the form just described. We develop in this paper yet another, global this time, optimality property for such diffusions.


\subsection{The probabilistic setting} \label{t.p.s}


In \hyperref[relative.entropy]{(\ref*{relative.entropy})} and throughout this paper, we are denoting by $\mathds{P}$ the unique probability measure on the space $\Omega = C([0,\infty);\mathds{R}^{n})$ of continuous, $\mathds{R}^{n}$-valued functions, under which the canonical coordinate process $X(t,\omega) = \omega(t)$, $t \geqslant 0$ has the property that
\begin{equation} \label{p.b.m}
W(t) \vcentcolon = X(t) - X(0) + \int_{0}^{t} \nabla \Psi\big( X(\theta) \big) \, \mathrm{d} \theta, \qquad t \geqslant 0
\end{equation}
is standard $\mathds{R}^{n}$-valued Brownian motion, and independent of the random variable $X(0)$ with distribution
\begin{equation}
\mathds{P}\big[ X(0) \in A \big] = \int_{A} p_{0}(x) \, \mathrm{d}x, \qquad A \in \mathscr{B}(\mathds{R}^{n}).
\end{equation}
The $\mathds{P}$-Brownian motion $(W(t))_{t \geqslant 0}$ of \hyperref[p.b.m]{(\ref*{p.b.m})} is adapted to, in fact generates, the canonical filtration $\mathds{F} = (\mathcal{F}(t))_{t \geqslant 0}$ with 
\begin{equation} \label{can.forw.filt}
\mathcal{F}(t) \vcentcolon = \boldsymbol{\sigma}\big( X(s) \colon \, 0 \leqslant s \leqslant t\big).
\end{equation}
By analogy to the terminology ``Wiener measure'', we call $\mathds{P}$ the ``Langevin--Smoluchowski measure'' associated with the potential $\Psi$.


\section{Reversal of time} \label{sec.three}


The densities $p(t, \, \cdot \,)$ and $q(\, \cdot \,)$ satisfy the forward Kolmogorov equations \hyperref[fp.fk.pde]{(\ref*{fp.fk.pde})} and \hyperref[fp.fk.pde.st]{(\ref*{fp.fk.pde.st})}, respectively. Whereas, their likelihood ratio $\ell(t, \, \cdot \,)$ in \hyperref[l.r.f]{(\ref*{l.r.f})} satisfies the \textit{backward} Kolmogorov equation \hyperref[b.k.pde]{(\ref*{b.k.pde})}. This suggests that, in the study of relative entropy and of its temporal dissipation, it might make sense to look at the underlying Langevin--Smoluchowski diffusion under time-reversal. Such an approach proved very fruitful in \cite{DPP89}, \cite{Pav89}, \cite{FJ16} and \cite{KST20a}; it will be important in our context here as well.

\smallskip

Thus, we fix an arbitrary terminal time $T \in (0,\infty)$ and consider the time-reversed process
\begin{equation} \label{can.pro}
\overline{X}(s) \vcentcolon = X(T-s) \, , \quad \qquad \overline{\mathcal{G}}(s) \vcentcolon = \boldsymbol{\sigma}\big( \overline{X}(u) \colon \, 0 \leqslant u \leqslant s \big) \, ; \quad \qquad 0 \leqslant s \leqslant T
\end{equation}
along with the filtration $\overline{\mathds{G}} = (\overline{\mathcal{G}}(s))_{0 \leqslant s \leqslant T}$ this process generates. Then standard theory on time-reversal shows that
\begin{equation} \label{bw.p.bm}
\overline{W}(s) \vcentcolon = \overline{X}(s) - \overline{X}(0) + \int_{0}^{s} \Big( \nabla \Psi \big( \overline{X}(u) \big) - \nabla L\big(T-u,\overline{X}(u)\big)\Big) \, \mathrm{d}u, \qquad 0 \leqslant s \leqslant T
\end{equation}
is, under $\mathds{P}$, a $\overline{\mathds{G}}$-standard Brownian motion with values in $\mathds{R}^{n}$ and independent of $\overline{X}(0) = X(T)$ (see, for instance, \cite{Foe85,Foe86}, \cite{HP86}, \cite{Mey94}, \cite{Nel01}, and \cite{Par86} for the classical results; an extensive presentation of the relevant facts regarding the time reversal of diffusion processes can be found in Appendix G of \cite{KST20b}). Here
\begin{equation} \label{log.lrfct}
L(t,x) \vcentcolon = \log \ell(t,x), \qquad (t,x) \in (0,\infty) \times \mathds{R}^{n}
\end{equation}
is the logarithm of the likelihood ratio function in \hyperref[l.r.f]{(\ref*{l.r.f})}; and on the strength of \hyperref[b.k.pde]{(\ref*{b.k.pde})}, this function solves the semilinear Schr\"odinger-type equation
\begin{equation} \label{sl.st.pde}
\partial L(t,x) = \tfrac{1}{2} \Delta L(t,x) - \big\langle \nabla L(t,x) \, , \nabla \Psi(x) \big\rangle + \tfrac{1}{2}  \vert \nabla L(t,x) \vert^{2}. 
\end{equation}

\smallskip

Another way to express this, is by saying that the so-called \textit{Hopf--Cole transform} $\ell = \mathrm{e}^{L}$ turns the semilinear equation \hyperref[sl.st.pde]{(\ref*{sl.st.pde})}, into the linear backward Kolmogorov equation \hyperref[b.k.pde]{(\ref*{b.k.pde})}. This observation is not new; it has been used in stochastic control to good effect by Fleming \cite{Fle77,Fle82}, Holland \cite{Hol77}, and in a context closer in spirit to this paper by Dai Pra and Pavon \cite{DPP90}, Dai Pra \cite{DP91}.


\section{A stochastic control problem} \label{sec.ascp.f}


Yet another way to cast the equation \hyperref[sl.st.pde]{(\ref*{sl.st.pde})}, is in the \textit{Hamilton--Jacobi--Bellman} form
\begin{equation} \label{sl.st.pde.m}
\partial L(t,x) =  \tfrac{1}{2} \Delta L(t,x) - \big\langle \nabla L(t,x) \, , \nabla \Psi(x) \big\rangle - \min_{g \in \mathds{R}^{n}} \Big( \big\langle \nabla L(t,x) \, , \, g \big\rangle + \tfrac{1}{2}  \vert g \vert^{2} \Big),
\end{equation}
where the minimization is attained by the gradient $g_{\ast} = - \nabla L(t,x)$. This, in turn, suggests a \textit{stochastic control problem} related to the backwards diffusive dynamics 
\begin{equation} \label{bwd.ls}
\mathrm{d} \overline{X}(s) = \Big( \nabla L\big(T-s,\overline{X}(s)\big) - \nabla \Psi \big( \overline{X}(s) \big) \Big) \, \mathrm{d}s + \mathrm{d}\overline{W}(s)
\end{equation}
of \hyperref[bw.p.bm]{(\ref*{bw.p.bm})}, which we formulate now as follows. 

\smallskip

For any measurable process $[0,T] \times \Omega \ni (t,\omega) \mapsto \boldsymbol{\gamma}(t,\omega) \in \mathds{R}^{n}$ such that the time-reversed process $\boldsymbol{\gamma}(T-s)$, $0 \leqslant s \leqslant T$ is adapted to the backward filtration $\overline{\mathds{G}}$ of \hyperref[can.pro]{(\ref*{can.pro})}, and which satisfies the condition
\begin{equation} \label{adap.pro}
\mathds{P}\bigg[ \int_{0}^{T} \vert \boldsymbol{\gamma}(T-s)\vert^{2} \, \mathrm{d}s < \infty \bigg] = 1,
\end{equation}
we consider the exponential $(\overline{\mathds{G}},\mathds{P})$-local martingale
\begin{equation} \label{exponential.martingale}
Z^{\boldsymbol{\gamma}}(s) \vcentcolon = \exp \bigg( \int_{0}^{s} \big\langle \boldsymbol{\gamma}(T-u) \, , \, \mathrm{d}\overline{W}(u) \big\rangle - \tfrac{1}{2} \int_{0}^{s} \vert \boldsymbol{\gamma}(T-u)\vert^{2} \, \mathrm{d}u \bigg) 
\, , \qquad 0 \leqslant s \leqslant T.
\end{equation}
We denote by $\boldsymbol{\Gamma}$ the collection of all processes $\boldsymbol{\gamma}$ as above, for which $Z^{\boldsymbol{\gamma}}$ is a true $(\overline{\mathds{G}},\mathds{P})$-martingale. This collection is not empty: it contains all such uniformly bounded processes $\boldsymbol{\gamma}$, and quite a few more (e.g., conditions of Novikov \cite[Corollary 3.5.13]{KS88} and Kazamaki \cite[Proposition VIII.1.14]{RY99}).

\smallskip

Now, for every $\boldsymbol{\gamma} \in \boldsymbol{\Gamma}$, we introduce an equivalent probability measure $\mathds{P}^{\boldsymbol{\gamma}} \sim \mathds{P}$ on path space, via
\begin{equation} \label{girs.dens.}
\frac{\mathrm{d}\mathds{P}^{\boldsymbol{\gamma}}}{\mathrm{d}\mathds{P}} \bigg\vert_{\overline{\mathcal{G}}(s)} 
= Z^{\boldsymbol{\gamma}}(s)
\, , \qquad 0 \leqslant s \leqslant T.
\end{equation}
Then, by the Girsanov theorem \cite[Theorem 3.5.1]{KS88}, the process
\begin{equation} \label{ti.re.ga.bm}
\overline{W}^{\boldsymbol{\gamma}}(s) \vcentcolon = \overline{W}(s) - \int_{0}^{s} \boldsymbol{\gamma}(T-u) \, \mathrm{d}u, \qquad 0 \leqslant s \leqslant T
\end{equation}
is standard $\mathds{R}^{n}$-valued $\mathds{P}^{\boldsymbol{\gamma}}$-Brownian motion of the filtration $\overline{\mathds{G}}$, thus independent of the random variable $\overline{X}(0) = X(T)$. Under the probability measure $\mathds{P}^{\boldsymbol{\gamma}}$, the backwards dynamics of \hyperref[bwd.ls]{(\ref*{bwd.ls})} take the form   
\begin{equation} \label{bwd.ls.sec}
\mathrm{d} \overline{X}(s) = \Big( \nabla L\big(T-s,\overline{X}(s)\big)  + \boldsymbol{\gamma}(T-s) - \nabla \Psi\big(\overline{X}(s)\big) \Big) \, \mathrm{d}s + \mathrm{d}\overline{W}^{\boldsymbol{\gamma}}(s);
\end{equation}
and it follows readily from this decomposition and the semilinear parabolic equation \hyperref[sl.st.pde]{(\ref*{sl.st.pde})}, that the process
\begin{equation} \label{pr.la}
M^{\boldsymbol{\gamma}}(s) \vcentcolon = L\big(T-s,\overline{X}(s)\big) + \tfrac{1}{2} \int_{0}^{s} \vert \boldsymbol{\gamma}(T-u) \vert^{2} \, \mathrm{d}u, \qquad 0 \leqslant s \leqslant T
\end{equation}
is a local $\overline{\mathds{G}}$-submartingale under $\mathds{P}^{\boldsymbol{\gamma}}$, with decomposition
\begin{equation} \label{gam.dr.te.va}
\mathrm{d}M^{\boldsymbol{\gamma}}(s) 
= \tfrac{1}{2} \big\vert \nabla L\big( T-s,\overline{X}(s)\big) + \boldsymbol{\gamma}(T-s) \big\vert^{2}  \, \mathrm{d}s + \Big \langle \nabla L\big(T-s,\overline{X}(s)\big) \, , \, \mathrm{d}\overline{W}^{\boldsymbol{\gamma}}(s) \Big\rangle.
\end{equation}
In fact, introducing the sequence
\begin{equation} \label{stop.time}
\sigma_{n} \vcentcolon = \inf \bigg\{ s \geqslant 0 \colon \, \int_{0}^{s} \Big( \big\vert \nabla L\big( T-u,\overline{X}(u)\big) \big\vert^{2}  + \vert \boldsymbol{\gamma}( T-u)\vert^{2} \Big) \, \mathrm{d}u \, \geqslant n \,  \bigg\} \wedge T
\, , \qquad n \in \mathds{N}_{0}
\end{equation}
of $\overline{\mathds{G}}$-stopping times with $\sigma_{n} \uparrow T$, we see that the stopped process $M^{\boldsymbol{\gamma}}(\, \cdot \, \wedge \sigma_{n})$ is a $\overline{\mathds{G}}$-submartingale under $\mathds{P}^{\boldsymbol{\gamma}}$, for every $n \in \mathds{N}_{0}$. In particular, we observe
\begingroup
\addtolength{\jot}{0.7em}
\begin{equation} \label{tr.st.subm.}
\begin{aligned} 
H\big( P(T) \, \vert \, \mathrm{Q} \big) &= \mathds{E}_{\mathds{P}}\big[ L\big(T,X(T)\big) \big] = \mathds{E}_{\mathds{P}^{\boldsymbol{\gamma}}}\big[ L\big(T,\overline{X}(0)\big) \big] \\
&\leqslant \mathds{E}_{\mathds{P}^{\boldsymbol{\gamma}}}\bigg[ L\big(T-\sigma_{n},\overline{X}(\sigma_{n})\big) + \tfrac{1}{2} \int_{0}^{\sigma_{n}} \vert \boldsymbol{\gamma}(T-u) \vert^{2} \, \mathrm{d}u \bigg],
\end{aligned}
\end{equation}
\endgroup
since $\mathds{P}^{\boldsymbol{\gamma}} = \mathds{P}$ on $\overline{\mathcal{G}}(0) = \boldsymbol{\sigma}(\overline{X}(0)) = \boldsymbol{\sigma}(X(T))$. Now \hyperref[tr.st.subm.]{(\ref*{tr.st.subm.})} holds for every $n \in \mathds{N}_{0}$, thus
\begin{equation} \label{tr.st.subm.lim}
H\big( P(T) \, \vert \, \mathrm{Q} \big) 
\leqslant \liminf_{n \rightarrow \infty} \ \mathds{E}_{\mathds{P}^{\boldsymbol{\gamma}}}\bigg[ L\big(T-\sigma_{n},\overline{X}(\sigma_{n})\big) + \tfrac{1}{2} \int_{0}^{\sigma_{n}} \vert \boldsymbol{\gamma}(T-u) \vert^{2} \, \mathrm{d}u \bigg].
\end{equation}

\smallskip

But as we remarked already, the minimum in \hyperref[sl.st.pde.m]{(\ref*{sl.st.pde.m})} is attained by $g_{\ast} = - \nabla L(t,x)$; likewise, the drift term in \hyperref[gam.dr.te.va]{(\ref*{gam.dr.te.va})} vanishes, if we select the process $\boldsymbol{\gamma}_{\ast} \in \boldsymbol{\Gamma}$ via 
\begin{equation} \label{score.pro}
\boldsymbol{\gamma}_{\ast}(t,\omega) \vcentcolon = - \nabla L\big(t,\omega(t)\big), \qquad \textnormal{ thus } \qquad \boldsymbol{\gamma}_{\ast}(T-s)  = - \nabla L\big(T-s,\overline{X}(s)\big), \qquad 0 \leqslant s \leqslant T.
\end{equation}
With this choice, the backwards dynamics of \hyperref[bwd.ls.sec]{(\ref*{bwd.ls.sec})} take the form
\begin{equation} \label{bwd.ls.sec.thir}
\mathrm{d} \overline{X}(s) =  - \nabla \Psi\big(\overline{X}(s)\big)  \, \mathrm{d}s + \mathrm{d}\overline{W}^{\boldsymbol{\gamma}_{\ast}}(s);
\end{equation}
that is, \textit{precisely of the Langevin--Smoluchowski type} \hyperref[ls.dp.sde]{(\ref*{ls.dp.sde})}, but now with the ``initial condition'' $\overline{X}(0) = X(T)$ and independent driving $\overline{\mathds{G}}$-Brownian motion $\overline{W}^{\boldsymbol{\gamma}_{\ast}}$, under $\mathds{P}^{\boldsymbol{\gamma}_{\ast}}$. Since $\mathds{P}^{\boldsymbol{\gamma}_{\ast}} = \mathds{P}$ holds on the $\boldsymbol{\sigma}$-algebra $\overline{\mathcal{G}}(0) = \boldsymbol{\sigma}(\overline{X}(0)) = \boldsymbol{\sigma}(X(T))$, the initial distribution of $\overline{X}(0)$ under $\mathds{P}^{\boldsymbol{\gamma}_{\ast}}$ is equal to $P(T)$. Furthermore, with $\boldsymbol{\gamma} = \boldsymbol{\gamma}_{\ast}$, the process of \hyperref[pr.la]{(\ref*{pr.la})}, \hyperref[gam.dr.te.va]{(\ref*{gam.dr.te.va})} becomes a $\mathds{P}^{\boldsymbol{\gamma}_{\ast}}$-local martingale, namely
\begin{equation} \label{loc.mart.}
M^{\boldsymbol{\gamma}_{\ast}}(s) 
= L\big(T,X(T)\big) + \int_{0}^{s} \Big \langle \nabla L\big(T-u,\overline{X}(u)\big) \, , \, \mathrm{d}\overline{W}^{\boldsymbol{\gamma}_{\ast}}(u) \Big\rangle, \qquad 0 \leqslant s \leqslant T;
\end{equation}
and we have equality in \hyperref[tr.st.subm.]{(\ref*{tr.st.subm.})}, thus also
\begin{equation} \label{tr.st.subm.lim.sec}
H\big( P(T) \, \vert \, \mathrm{Q} \big) 
= \lim_{n \rightarrow \infty} \ \mathds{E}_{\mathds{P}^{\boldsymbol{\gamma}_{\ast}}}\bigg[ L\big(T-\sigma_{n},\overline{X}(\sigma_{n})\big) + \tfrac{1}{2} \int_{0}^{\sigma_{n}} \vert \boldsymbol{\gamma}_{\ast}(T-u) \vert^{2} \, \mathrm{d}u \bigg].
\end{equation}

\smallskip

We conclude that the infimum over $\boldsymbol{\gamma} \in \boldsymbol{\Gamma}$ of the right-hand side in \hyperref[tr.st.subm.lim]{(\ref*{tr.st.subm.lim})} is attained by the process $\boldsymbol{\gamma}_{\ast}$ of \hyperref[score.pro]{(\ref*{score.pro})}, which gives rise to the Langevin--Smoluchowski dynamics \hyperref[bwd.ls.sec.thir]{(\ref*{bwd.ls.sec.thir})} for the time-reversed process $\overline{X}(s) = X(T-s)$, $0 \leqslant s \leqslant T$, under $\mathds{P}^{\boldsymbol{\gamma}_{\ast}}$. We formalize this discussion as follows.

\begin{theorem} \label{m.theo.one} Consider the stochastic control problem of minimizing over the class $\boldsymbol{\Gamma}$ of measurable, adapted processes $\boldsymbol{\gamma}$ satisfying \textnormal{\hyperref[adap.pro]{(\ref*{adap.pro})}} and inducing an exponential martingale $Z^{\boldsymbol{\gamma}}$ in \textnormal{\hyperref[exponential.martingale]{(\ref*{exponential.martingale})}}, with the notation of \textnormal{\hyperref[stop.time]{(\ref*{stop.time})}} and with the backwards dynamics of \textnormal{\hyperref[bwd.ls.sec]{(\ref*{bwd.ls.sec})}}, the expected cost 
\begin{equation} \label{expected.cost}
\mathcal{I}(\boldsymbol{\gamma}) \vcentcolon = \liminf_{n \rightarrow \infty} \ \mathds{E}_{\mathds{P}^{\boldsymbol{\gamma}}}\bigg[ L\big(T-\sigma_{n},\overline{X}(\sigma_{n})\big) + \tfrac{1}{2} \int_{0}^{\sigma_{n}} \vert \boldsymbol{\gamma}(T-u) \vert^{2} \, \mathrm{d}u \bigg].
\end{equation}

\smallskip

Under the assumptions of \textnormal{\hyperref[the.setting]{Section \ref*{the.setting}}}, the infimum $\inf_{\boldsymbol{\gamma} \in \boldsymbol{\Gamma}} \mathcal{I}(\boldsymbol{\gamma})$ is equal to the relative entropy $H( P(T) \, \vert \, \mathrm{Q})$ and is attained by the ``score process'' $\boldsymbol{\gamma}_{\ast}$ of \textnormal{\hyperref[score.pro]{(\ref*{score.pro})}}. This choice leads to the backwards Langevin--Smoluchowski dynamics \textnormal{\hyperref[bwd.ls.sec.thir]{(\ref*{bwd.ls.sec.thir})}}, and with $\boldsymbol{\gamma} = \boldsymbol{\gamma}_{\ast}$ the limit in \textnormal{\hyperref[expected.cost]{(\ref*{expected.cost})}} exists as in \textnormal{\hyperref[tr.st.subm.lim.sec]{(\ref*{tr.st.subm.lim.sec})}}. 
\begin{proof} It only remains to check that the minimizing process of \hyperref[score.pro]{(\ref*{score.pro})} belongs indeed to the collection $\boldsymbol{\Gamma}$ of admissible processes. By its definition, this process $\boldsymbol{\gamma}_{\ast}$ is measurable, and its time-reversal is adapted to the backward filtration $\overline{\mathds{G}}$ of \hyperref[can.pro]{(\ref*{can.pro})}. Theorem 4.1 in \cite{KST20a} gives
\begin{equation} \label{finite.kst}
\mathds{E}_{\mathds{P}}\bigg[ \int_{0}^{T} \big\vert \nabla L\big(T-u,\overline{X}(u)\big) \big\vert^{2} \, \mathrm{d}u  \bigg]
= \mathds{E}_{\mathds{P}}\bigg[ \int_{0}^{T} \big\vert \nabla L\big(\theta,X(\theta)\big) \big\vert^{2} \, \mathrm{d}\theta  \bigg] < \infty,
\end{equation}
which implies \textit{a fortiori} that the condition in \hyperref[adap.pro]{(\ref*{adap.pro})} is satisfied for $\boldsymbol{\gamma} = \boldsymbol{\gamma}_{\ast}$.

\smallskip

We must also show that the process $Z^{\boldsymbol{\gamma}_{\ast}}$ defined in the manner of \hyperref[exponential.martingale]{(\ref*{exponential.martingale})}, is a true martingale. A very mild dose of stochastic calculus leads to
\begin{equation}
\mathrm{d} L\big(T-s,\overline{X}(s)\big) 
= \Big\langle \nabla L\big(T-s,\overline{X}(s)\big) \, , \, \mathrm{d}\overline{W}(s) \Big\rangle  + \tfrac{1}{2} \big\vert\nabla L\big(T-s,\overline{X}(s)\big)\big\vert^{2} \, \mathrm{d}s 
\end{equation}
on account of \hyperref[sl.st.pde]{(\ref*{sl.st.pde})}, \hyperref[bwd.ls]{(\ref*{bwd.ls})}. Therefore, we have
\begingroup
\addtolength{\jot}{0.7em}
\begin{align}
&\int_{0}^{s} \big\langle \boldsymbol{\gamma}_{\ast}(T-u) \, , \, \mathrm{d}\overline{W}(u) \big\rangle - \tfrac{1}{2} \int_{0}^{s} \vert \boldsymbol{\gamma}_{\ast}(T-u)\vert^{2} \, \mathrm{d}u \\
& \qquad = - \int_{0}^{s} \Big\langle \nabla L\big(T-u,\overline{X}(u)\big) \, , \, \mathrm{d}\overline{W}(u) \Big\rangle - \tfrac{1}{2} \int_{0}^{s} \big\vert \nabla L\big(T-u,\overline{X}(u)\big) \big\vert^{2} \, \mathrm{d}u \\
& \qquad = L\big(T,X(T)\big) - L\big(T-s,X(T-s)\big) = \log \bigg( \frac{\ell\big(T,X(T)\big)}{\ell\big(T-s,X(T-s)\big)} \bigg),
\end{align}
\endgroup
which expresses the exponential process of \hyperref[exponential.martingale]{(\ref*{exponential.martingale})} with $\boldsymbol{\gamma} = \boldsymbol{\gamma}_{\ast}$ as
\begin{equation}
Z^{\boldsymbol{\gamma}_{\ast}}(s) = \frac{\ell\big(T,X(T)\big)}{\ell\big(T-s,X(T-s)\big)} \, , \qquad 0 \leqslant s \leqslant T.
\end{equation}

\smallskip

Now, let us argue that the process $Z^{\boldsymbol{\gamma}_{\ast}}$ is a true $(\overline{\mathds{G}},\mathds{P})$-martingale. It is a positive local martingale, thus a supermartingale. It will be a martingale, if it has constant expectation. But $Z^{\boldsymbol{\gamma}_{\ast}}(0) \equiv 1$, so it is enough to show that $\mathds{E}_{\mathds{P}}[Z^{\boldsymbol{\gamma}_{\ast}}(T)] = 1$. Let us denote by $P(s,y;t,\xi)$ the transition kernel of the Langevin--Smoluchowski dynamics, so that $\mathds{P}[X(s) \in \mathrm{d}y, X(t) \in \mathrm{d}\xi] = p(s,y) \, P(s,y;t,\xi) \, \mathrm{d}y \, \mathrm{d}\xi$ for $0 \leqslant s < t \leqslant T$ and $(y,\xi) \in \mathds{R}^{n} \times \mathds{R}^{n}$. Then the invariance of $\mathrm{Q}$ gives
\begin{equation} \label{invar.gam}
\int_{\mathds{R}^{n}} q(y) \,  P(0,y;T,\xi) \, \mathrm{d}y = q(\xi) \, , \qquad \xi \in \mathds{R}^{n};
\end{equation}
consequently
\begingroup
\addtolength{\jot}{0.7em}
\begin{align}
\mathds{E}_{\mathds{P}}\big[Z^{\boldsymbol{\gamma}_{\ast}}(T)\big] 
&= \mathds{E}_{\mathds{P}}\bigg[ \frac{p\big(T,X(T)\big)}{q\big(X(T)\big)} \frac{q\big(X(0)\big)}{p\big(0,X(0)\big)} \bigg] \label{invar.gam.a} \\
&= \int_{\mathds{R}^{n}} \int_{\mathds{R}^{n}} \frac{p(T,\xi)}{q(\xi)} \frac{q(y)}{p(0,y)} \, p(0,y) \, P(0,y;T,\xi) \, \mathrm{d}y \, \mathrm{d}\xi \label{invar.gam.b} \\
&= \int_{\mathds{R}^{n}}  \frac{p(T,\xi)}{q(\xi)}  \bigg( \int_{\mathds{R}^{n}} q(y) \, P(0,y;T,\xi) \, \mathrm{d}y \bigg) \, \mathrm{d}\xi \label{invar.gam.c} \\
&= \int_{\mathds{R}^{n}} p(T,\xi) \, \mathrm{d}\xi = 1, \label{invar.gam.d}
\end{align}
\endgroup
implying that $Z^{\boldsymbol{\gamma}_{\ast}}$ is a true martingale and completing the proof of \hyperref[m.theo.one]{Theorem \ref*{m.theo.one}}.
\end{proof}
\end{theorem}

Results related to \hyperref[m.theo.one]{Theorem \ref*{m.theo.one}} have been established in \cite{DPP89,DPP90,FS85,Pav89,PW91}.

\begin{remark}[\textsc{Reincarnation of time-marginals}] Let us denote by $\overline{P}_{\ast}(s)$ the distribution of the random variable $\overline{X}(s) = X(T-s)$ under the probability measure $\mathds{P}^{\boldsymbol{\gamma}_{\ast}}$, for $0 \leqslant s \leqslant T$. Since $(\overline{X}(s))_{0 \leqslant s \leqslant T}$ is under $\mathds{P}^{\boldsymbol{\gamma}_{\ast}}$ a Langevin--Smoluchowski diffusion in its own right, we deduce
\begin{equation} \label{uniq.l.s.d}
\overline{P}_{\ast}(s) = P(T+s), \qquad 0 \leqslant s \leqslant T
\end{equation}
on the strength of uniqueness in distribution for the Langevin--Smoluchowski flow, and of its time-homogeneity. In other words, the branch $(P(T+s))_{0 \leqslant s \leqslant T}$ of the original Langevin--Smoluchowski curve of time-marginals, gets ``reincarnated'' as $(\overline{P}_{\ast}(s))_{0 \leqslant s \leqslant T}$, the curve of time-marginals arising from the solution of the stochastic control problem in \hyperref[m.theo.one]{Theorem \ref*{m.theo.one}}. But now, under the probability measure $\mathds{P}^{\boldsymbol{\gamma}_{\ast}}$, the states of the Langevin--Smoluchowski diffusion $(\overline{X}(s))_{0 \leqslant s \leqslant T}$ corresponding to the curve $(\overline{P}_{\ast}(s))_{0 \leqslant s \leqslant T}$ traverse the time interval $[0,T]$ in the opposite temporal direction. 
\end{remark}


\subsection{Discussion: entropic interpretation of the expected cost when \texorpdfstring{$\mathrm{Q}(\mathds{R}^{n}) < \infty$}{Q is a probability measure}} \label{int.one}


Let us observe from \hyperref[exponential.martingale]{(\ref*{exponential.martingale})} -- \hyperref[ti.re.ga.bm]{(\ref*{ti.re.ga.bm})} that
\begin{equation} \label{log.girs.dens.rel.ent.}
\log \bigg( \frac{\mathrm{d}\mathds{P}^{\boldsymbol{\gamma}}}{\mathrm{d}\mathds{P}} \bigg\vert_{\overline{\mathcal{G}}(\sigma_{n})} \bigg)
=  \int_{0}^{\sigma_{n}} \big\langle \boldsymbol{\gamma}(T-u) \, , \, \mathrm{d}\overline{W}^{\boldsymbol{\gamma}}(u) \big\rangle + \tfrac{1}{2} \int_{0}^{\sigma_{n}} \vert \boldsymbol{\gamma}(T-u)\vert^{2} \, \mathrm{d}u 
\end{equation}
holds for every $\boldsymbol{\gamma} \in \boldsymbol{\Gamma}$ and $n \in \mathds{N}_{0}$. Thus, as the $\mathds{P}^{\boldsymbol{\gamma}}$-expectation of the stochastic integral in \hyperref[log.girs.dens.rel.ent.]{(\ref*{log.girs.dens.rel.ent.})} vanishes, the expected quadratic cost, or ``energy'', term in \hyperref[expected.cost]{(\ref*{expected.cost})} is itself a relative entropy:
\begin{equation}
\mathds{E}_{\mathds{P}^{\boldsymbol{\gamma}}}\bigg[ \tfrac{1}{2} \int_{0}^{\sigma_{n}} \vert \boldsymbol{\gamma}(T-u)\vert^{2} \, \mathrm{d}u \bigg]
= \mathds{E}_{\mathds{P}^{\boldsymbol{\gamma}}}\bigg[ \log \bigg( \frac{\mathrm{d}\mathds{P}^{\boldsymbol{\gamma}}}{\mathrm{d}\mathds{P}} \bigg\vert_{\overline{\mathcal{G}}(\sigma_{n})} \bigg)\bigg].
\end{equation}

\smallskip

By contrast, when $\mathrm{Q}$ is a probability measure on $\mathscr{B}(\mathds{R}^{n})$, and denoting by $\mathds{Q}$ the probability measure induced on $\Omega = C([0,\infty);\mathds{R}^{n})$ by the canonical process driven by the dynamics \hyperref[ls.dp.sde]{(\ref*{ls.dp.sde})} with $\mathrm{Q}$ as the distribution of $X(0)$, the first term in \hyperref[expected.cost]{(\ref*{expected.cost})} can be cast as 
\begingroup
\addtolength{\jot}{0.7em}
\begin{align}
\mathds{E}_{\mathds{P}^{\boldsymbol{\gamma}}}\big[ L\big(T-\sigma_{n},\overline{X}(\sigma_{n})\big) \big] 
&= \mathds{E}_{\mathds{P}^{\boldsymbol{\gamma}}}\bigg[ \log \bigg( \frac{\mathrm{d}\mathds{P}}{\mathrm{d}\mathds{Q}} \bigg\vert_{\boldsymbol{\sigma}(\overline{X}(\sigma_{n}))} \bigg)\bigg] \label{narrent.a} \\
&= \mathds{E}_{\mathds{P}^{\boldsymbol{\gamma}}}\bigg[ \log \bigg( \frac{\mathrm{d}\mathds{P}^{\boldsymbol{\gamma}}}{\mathrm{d}\mathds{Q}} \bigg\vert_{\boldsymbol{\sigma}(\overline{X}(\sigma_{n}))} \bigg)\bigg] 
- \mathds{E}_{\mathds{P}^{\boldsymbol{\gamma}}}\bigg[ \log \bigg( \frac{\mathrm{d}\mathds{P}^{\boldsymbol{\gamma}}}{\mathrm{d}\mathds{P}} \bigg\vert_{\boldsymbol{\sigma}(\overline{X}(\sigma_{n}))} \bigg)\bigg]. \label{narrent.b}
\end{align}
\endgroup
It follows that, in this case, the expected cost of \hyperref[expected.cost]{(\ref*{expected.cost})} is equal to the sum $H_{n}^{\boldsymbol{\gamma}} + D_{n}^{\boldsymbol{\gamma}}$ of two non-negative quantities:
\begin{equation}
H_{n}^{\boldsymbol{\gamma}} \vcentcolon = \mathds{E}_{\mathds{P}^{\boldsymbol{\gamma}}}\bigg[ \log \bigg( \frac{\mathrm{d}\mathds{P}^{\boldsymbol{\gamma}}}{\mathrm{d}\mathds{Q}} \bigg\vert_{\boldsymbol{\sigma}(\overline{X}(\sigma_{n}))} \bigg)\bigg],
\end{equation}
the relative entropy of the probability measure $\mathds{P}^{\boldsymbol{\gamma}}$ with respect to the probability measure $\mathds{Q}$ when both are restricted to the $\boldsymbol{\sigma}$-algebra generated by the random variable $\overline{X}(\sigma_{n})$; and
\begin{equation} \label{ent.diff}
D_{n}^{\boldsymbol{\gamma}} \vcentcolon = \mathds{E}_{\mathds{P}^{\boldsymbol{\gamma}}}\bigg[ \log \bigg( \frac{\mathrm{d}\mathds{P}^{\boldsymbol{\gamma}}}{\mathrm{d}\mathds{P}} \bigg\vert_{\overline{\mathcal{G}}(\sigma_{n})} \bigg)\bigg] 
- \mathds{E}_{\mathds{P}^{\boldsymbol{\gamma}}}\bigg[ \log \bigg( \frac{\mathrm{d}\mathds{P}^{\boldsymbol{\gamma}}}{\mathrm{d}\mathds{P}} \bigg\vert_{\boldsymbol{\sigma}(\overline{X}(\sigma_{n}))} \bigg)\bigg],
\end{equation}
the difference between the relative entropies of the probability measure $\mathds{P}^{\boldsymbol{\gamma}}$ with respect to the probability measure $\mathds{P}$, when restricted to the $\boldsymbol{\sigma}$-algebra generated by the collection of random variables $(\overline{X}(u \wedge \sigma_{n}))_{0 \leqslant u \leqslant T}$ and by the random variable $\overline{X}(\sigma_{n})$, respectively. The difference in \hyperref[ent.diff]{(\ref*{ent.diff})} is non-negative, because conditioning on a smaller $\boldsymbol{\sigma}$-algebra can only decrease the relative entropy; this difference can be thought of as an ``entropic cost of time-reversal''.

\smallskip

It develops from this discussion that the expected cost on the right-hand side of \hyperref[expected.cost]{(\ref*{expected.cost})} is non-negative, when $\mathrm{Q}$ is a probability measure.


\section{From local to square-integrable martingales} \label{sec.ltm}


Whenever the process $M^{\boldsymbol{\gamma}}$ of \hyperref[pr.la]{(\ref*{pr.la})}, \hyperref[gam.dr.te.va]{(\ref*{gam.dr.te.va})}  happens to be a true submartingale under $\mathds{P}^{\boldsymbol{\gamma}}$ (as, for instance, with $\boldsymbol{\gamma} \equiv 0$ on account of Theorem 4.1 in \cite{KST20a}), the expected cost \hyperref[expected.cost]{(\ref*{expected.cost})} takes the form
\begin{equation} 
\mathds{E}_{\mathds{P}^{\boldsymbol{\gamma}}}\bigg[ L\big(0,X(0)\big) + \tfrac{1}{2} \int_{0}^{T} \vert \boldsymbol{\gamma}(T-u) \vert^{2} \, \mathrm{d}u \bigg].
\end{equation}
Likewise, we derive from \hyperref[m.theo.one]{Theorem \ref*{m.theo.one}} the identity
\begin{equation} 
H\big( P(T) \, \vert \, \mathrm{Q} \big) 
= \mathds{E}_{\mathds{P}^{\boldsymbol{\gamma}_{\ast}}}\bigg[ L\big(0,X(0)\big) + \tfrac{1}{2} \int_{0}^{T} \big\vert \nabla L\big(T-u,\overline{X}(u)\big) \big\vert^{2} \, \mathrm{d}u \bigg],
\end{equation}
whenever the process $M^{\boldsymbol{\gamma}_{\ast}}$ of \hyperref[loc.mart.]{(\ref*{loc.mart.})} is a true $\mathds{P}^{\boldsymbol{\gamma}_{\ast}}$-martingale. This is the case, for instance, whenever
\begin{equation} \label{cond.true.mart}
\mathds{E}_{\mathds{P}^{\boldsymbol{\gamma}_{\ast}}} \Big[ \big\langle M^{\boldsymbol{\gamma}_{\ast}}, M^{\boldsymbol{\gamma}_{\ast}} \big\rangle(T) \Big] = \mathds{E}_{\mathds{P}^{\boldsymbol{\gamma}_{\ast}}}\bigg[ \int_{0}^{T} \big\vert \nabla L\big(T-u,\overline{X}(u)\big) \big\vert^{2} \, \mathrm{d}u  \bigg] < \infty
\end{equation}
holds; and then the stochastic integral in \hyperref[loc.mart.]{(\ref*{loc.mart.})} is an $L^{2}(\mathds{P}^{\boldsymbol{\gamma}_{\ast}})$-bounded martingale (see, for instance, \cite[Proposition 3.2.10]{KS88} or \cite[Corollary IV.1.25]{RY99}). Using \hyperref[uniq.l.s.d]{(\ref*{uniq.l.s.d})}, we can express the expectation of \hyperref[cond.true.mart]{(\ref*{cond.true.mart})} more explicitly as 
\begin{equation} \label{cond.true.mart.ex}
\int_{\mathds{R}^{n}} \int_{0}^{T} \vert \nabla \log \ell(t,x) \vert^{2} \, p(2T-t,x) \, \mathrm{d}t  \, \mathrm{d}x.
\end{equation}
The shift in the temporal variable makes it difficult to check whether the quantity in \hyperref[cond.true.mart.ex]{(\ref*{cond.true.mart.ex})} is finite. At least, we have not been able to apply directly arguments similar to those in Theorem 4.1 of \cite{KST20a}, where the expectation \hyperref[cond.true.mart]{(\ref*{cond.true.mart})} is taken with respect to the probability measure $\mathds{P}$, in the manner of \hyperref[finite.kst]{(\ref*{finite.kst})} (and thus, the argument $2T-t$ in \hyperref[cond.true.mart.ex]{(\ref*{cond.true.mart.ex})} is replaced by $t$). This problem is consonant with the fact that the term in \hyperref[narrent.a]{(\ref*{narrent.a})}, \hyperref[narrent.b]{(\ref*{narrent.b})} is not quite a relative entropy, but a linear combination of relative entropies.

\smallskip

The goal of this section is to find a square-integrable $\mathds{P}^{\boldsymbol{\gamma}_{\ast}}$-martingale $\overline{M}^{\boldsymbol{\gamma}_{\ast}}$, which is closely related to the local martingale $M^{\boldsymbol{\gamma}_{\ast}}$ of \hyperref[loc.mart.]{(\ref*{loc.mart.})}. The idea is to correct the shift in the temporal variable appearing in \hyperref[cond.true.mart.ex]{(\ref*{cond.true.mart.ex})}, by reversing time once again.

\smallskip

First, we need to introduce some notation. We denote by $\overline{p}_{\ast}(s,\, \cdot \,)$ the probability density function of the random variable $\overline{X}(s) = X(T-s)$ under the probability measure $\mathds{P}^{\boldsymbol{\gamma}_{\ast}}$, for $0 \leqslant s \leqslant T$. From \hyperref[uniq.l.s.d]{(\ref*{uniq.l.s.d})}, we deduce the relation
\begin{equation} \label{uniq.l.s.d.p}
\overline{p}_{\ast}(s,x) = p(T+s,x), \qquad (s,x) \in [0,T] \times \mathds{R}^{n}.
\end{equation}
The associated likelihood ratio function and its logarithm are defined respectively by
\begin{equation} \label{log.lrfct.n}
\lambda(s,x) \vcentcolon = \frac{\overline{p}_{\ast}(s,x)}{q(x)} \, , \qquad \Lambda(s,x) \vcentcolon = \log \lambda(s,x) = L(T+s,x) \, , \qquad (s,x) \in [0,T] \times \mathds{R}^{n}.
\end{equation}
From the definition \hyperref[log.lrfct.n]{(\ref*{log.lrfct.n})}, and the equations \hyperref[sl.st.pde]{(\ref*{sl.st.pde})}, \hyperref[sl.st.pde.m]{(\ref*{sl.st.pde.m})}, we see that the function $(s,x) \mapsto \Lambda(s,x)$ satisfies again the semilinear Schr\"odinger-type equation
\begin{equation} \label{sl.st.pde.m.sec.three}
- \partial \Lambda(s,x) + \tfrac{1}{2} \Delta \Lambda(s,x) - \big\langle \nabla \Lambda(s,x) \, , \nabla \Psi(x) \big\rangle = \min_{b \in \mathds{R}^{n}} \Big( \big\langle \nabla \Lambda(s,x) \, , \, b \big\rangle + \tfrac{1}{2}  \vert b \vert^{2} \Big).
\end{equation}

\smallskip

In the setting introduced above, the relative entropy with respect to $\mathrm{Q}$ of the distribution $\overline{P}_{\ast}(s)$ of $\overline{X}(s)$ under $\mathds{P}^{\boldsymbol{\gamma}_{\ast}}$ is
\begin{equation} \label{rel.ent.sec}
H\big( \overline{P}_{\ast}(s) \, \vert \, \mathrm{Q} \big) 
= \mathds{E}_{\mathds{P}^{\boldsymbol{\gamma}_{\ast}}}\big[  \Lambda\big(s,\overline{X}(s)\big) \big]
= \int_{\mathds{R}^{n}} \log \bigg(\frac{\overline{p}_{\ast}(s,x)}{q(x)}\bigg) \, \overline{p}_{\ast}(s,x) \, \mathrm{d}x, \qquad 0 \leqslant s \leqslant T.
\end{equation}
Again, the assumption that $H(P(0) \, \vert \, \mathrm{Q})$ is finite, and the decrease of the relative entropy function $[0,\infty) \ni t \mapsto H( P(t) \, \vert \, \mathrm{Q}) \in (-\infty,\infty]$, imply that the relative entropy in \hyperref[rel.ent.sec]{(\ref*{rel.ent.sec})} is finite for all $0 \leqslant s \leqslant T$.

\smallskip

Finally, the relative Fisher information of $\overline{P}_{\ast}(s)$ with respect to $\mathrm{Q}$ is defined as
\begin{equation}  \label{fish.inf}
I\big( \overline{P}_{\ast}(s) \, \vert \, \mathrm{Q}\big) \vcentcolon = \mathds{E}_{\mathds{P}^{\boldsymbol{\gamma}_{\ast}}}\Big[ \, \big\vert \nabla \Lambda\big(s,\overline{X}(s)\big) \big\vert^{2} \, \Big] = \int_{\mathds{R}^{n}} \big\vert \nabla \Lambda(s,x) \big\vert^{2} \, \overline{p}_{\ast}(s,x) \, \mathrm{d}x, \qquad 0 \leqslant s \leqslant T.
\end{equation}

\subsection{Reversing time once again} \label{rev.time.again}

Let us consider on the filtered probability space $(\Omega,\mathds{F},\mathds{P}^{\boldsymbol{\gamma}_{\ast}})$ the canonical process $(X(t))_{0 \leqslant t \leqslant T}$, whose time-reversal \hyperref[can.pro]{(\ref*{can.pro})} satisfies the backwards Langevin--Smoluchowski dynamics \hyperref[bwd.ls.sec.thir]{(\ref*{bwd.ls.sec.thir})}. Reversing time once again, we find that $X(t) = \overline{X}(T-t)$, $0 \leqslant t \leqslant T$ satisfies the stochastic differential equation 
\begin{equation} \label{forward.dyn}
\mathrm{d} X(t) = \Big( \nabla \Lambda\big(T-t,X(t)\big) - \nabla \Psi\big(X(t)\big) \Big) \, \mathrm{d}t + \mathrm{d}W^{\boldsymbol{\gamma}_{\ast}}(t), \qquad 0 \leqslant t \leqslant T,
\end{equation}
where the process
\begin{equation} 
W^{\boldsymbol{\gamma}_{\ast}}(t) \vcentcolon = \overline{W}^{\boldsymbol{\gamma}_{\ast}}(T-t) - \overline{W}^{\boldsymbol{\gamma}_{\ast}}(T) - \int_{0}^{t} \nabla \log \overline{p}_{\ast}\big(T-\theta,X(\theta)\big) \, \mathrm{d}\theta
, \qquad 0 \leqslant t \leqslant T
\end{equation}
is Brownian motion on $(\Omega,\mathds{F},\mathds{P}^{\boldsymbol{\gamma}_{\ast}})$. We recall here Proposition 4.1 from \cite{KST20b}. Comparing the equation \hyperref[forward.dyn]{(\ref*{forward.dyn})} with \hyperref[ls.dp.sde]{(\ref*{ls.dp.sde})}, we see that the $\mathds{P}^{\boldsymbol{\gamma}_{\ast}}$-Brownian motion $(W^{\boldsymbol{\gamma}_{\ast}}(t))_{0 \leqslant t \leqslant T}$ and the $\mathds{P}$-Brownian motion $(W(t))_{0 \leqslant t \leqslant T}$ are related via
\begin{equation} 
W(t) = W^{\boldsymbol{\gamma}_{\ast}}(t) + \int_{0}^{t} \nabla \Lambda\big(T-\theta,X(\theta)\big) \, \mathrm{d}\theta, \qquad 0 \leqslant t \leqslant T.
\end{equation}

\subsection{The dynamics of the relative entropy process}

We look now at the \textit{relative entropy process}
\begin{equation} \label{rel.ent.pro}
\Lambda\big(T-t,X(t)\big) = \log \bigg( \frac{\overline{p}_{\ast}\big(T-t,X(t)\big)}{q\big(X(t)\big)} \bigg), \qquad 0 \leqslant t \leqslant T
\end{equation}
on $(\Omega,\mathds{F},\mathds{P}^{\boldsymbol{\gamma}_{\ast}})$. Applying It\^{o}'s formula and using the equation \hyperref[sl.st.pde.m.sec.three]{(\ref*{sl.st.pde.m.sec.three})}, together with the forward dynamics \hyperref[forward.dyn]{(\ref*{forward.dyn})}, we obtain the following result.

\begin{proposition} \label{rel.ent.prop.} On the filtered probability space $(\Omega,\mathds{F},\mathds{P}^{\boldsymbol{\gamma}_{\ast}})$, the relative entropy process \textnormal{\hyperref[rel.ent.pro]{(\ref*{rel.ent.pro})}} is a submartingale with stochastic differential
\begin{equation}  \label{loc.m.c}
\mathrm{d}\Lambda\big(T-t,X(t)\big) 
= \tfrac{1}{2} \big\vert \nabla \Lambda\big( T-t,X(t)\big) \big\vert^{2}  \, \mathrm{d}t + \Big \langle \nabla \Lambda\big(T-t,X(t)\big) \, , \, \mathrm{d}W^{\boldsymbol{\gamma}_{\ast}}(t) \Big\rangle
\end{equation}
for $0 \leqslant t \leqslant T$. In particular, the process
\begin{equation} \label{loc.m.a}
\overline{M}^{\boldsymbol{\gamma}_{\ast}}(t) \vcentcolon = \Lambda\big(T-t,X(t)\big) - \Lambda\big(T,X(0)\big) - \tfrac{1}{2} \int_{0}^{t} \big\vert \nabla \Lambda\big( T-\theta,X(\theta)\big) \big\vert^{2} \, \mathrm{d}\theta, \qquad 0 \leqslant t \leqslant T    
\end{equation}
is an $L^{2}(\mathds{P}^{\boldsymbol{\gamma}_{\ast}})$-bounded martingale, with stochastic integral representation
\begin{equation} \label{loc.m.b}
\overline{M}^{\boldsymbol{\gamma}_{\ast}}(t) = \int_{0}^{t} \Big \langle \nabla \Lambda\big(T-\theta,X(\theta)\big) \, , \, \mathrm{d}W^{\boldsymbol{\gamma}_{\ast}}(\theta) \Big\rangle, \qquad 0 \leqslant t \leqslant T.    
\end{equation}
\begin{proof} The last thing we need to verify for the proof of \hyperref[rel.ent.prop.]{Proposition \ref*{rel.ent.prop.}}, is that
\begin{equation} \label{true.exp.a}
  \mathds{E}_{\mathds{P}^{\boldsymbol{\gamma}_{\ast}}} \Big[ \big\langle \overline{M}^{\boldsymbol{\gamma}_{\ast}}, \overline{M}^{\boldsymbol{\gamma}_{\ast}} \big\rangle(T) \Big] 
= \mathds{E}_{\mathds{P}^{\boldsymbol{\gamma}_{\ast}}} \bigg[ \int_{0}^{T} \big\vert \nabla \Lambda\big( T-t,X(t)\big) \big\vert^{2} \, \mathrm{d}t\bigg] < \infty.
\end{equation}
We observe that the expectation in \hyperref[true.exp.a]{(\ref*{true.exp.a})} is equal to
\begin{equation} \label{true.exp.b}
   \mathds{E}_{\mathds{P}^{\boldsymbol{\gamma}_{\ast}}}\bigg[ \int_{0}^{T} \big\vert \nabla \Lambda\big(s,\overline{X}(s)\big) \big\vert^{2} \, \mathrm{d}s\bigg] 
= \mathds{E}_{\mathds{P}}\bigg[ \int_{T}^{2T} \big\vert \nabla L\big(t,X(t)\big) \big\vert^{2} \, \mathrm{d}t\bigg].
\end{equation}
This is because \hyperref[can.pro]{(\ref*{can.pro})}, \hyperref[log.lrfct.n]{(\ref*{log.lrfct.n})} give $\nabla \Lambda(s,\overline{X}(s)) = \nabla L(t,X(2T-t))$ with $t = T+s \in [T,2T]$; and because the $\mathds{P}^{\boldsymbol{\gamma}_{\ast}}$-distribution of $X(2T-t) = \overline{X}(s)$ is the same as the $\mathds{P}$-distribution of $X(T+s) = X(t)$, on account of \hyperref[uniq.l.s.d]{(\ref*{uniq.l.s.d})}. But, as \hyperref[finite.kst]{(\ref*{finite.kst})} holds for any finite time horizon $T > 0$, the quantity in \hyperref[true.exp.b]{(\ref*{true.exp.b})} is finite as well.
\end{proof}
\end{proposition}

\subsection{Relative entropy dissipation}

Exploiting the trajectorial evolution of the relative entropy process \hyperref[rel.ent.pro]{(\ref*{rel.ent.pro})}, provided by \hyperref[rel.ent.prop.]{Proposition \ref*{rel.ent.prop.}}, allows us to derive some immediate consequences on the decrease of the relative entropy function $[0,T] \ni s \longmapsto H( \overline{P}_{\ast}(s) \, \vert \, \mathrm{Q}) \in (-\infty,\infty)$ and its rate of dissipation. The submartingale-property of the relative entropy process \hyperref[rel.ent.pro]{(\ref*{rel.ent.pro})} shows once more, that this function is non-decreasing. More precisely, we have the following rate of change for the relative entropy.

\begin{corollary} \label{cor.a} For all $s,s_{0} \geqslant 0$, we have 
\begin{equation} \label{cor.a.b}
H\big( \overline{P}_{\ast}(s) \, \vert \, \mathrm{Q} \big) - H\big( \overline{P}_{\ast}(s_{0}) \, \vert \, \mathrm{Q} \big) 
= - \tfrac{1}{2} \int_{s_{0}}^{s}  I\big( \overline{P}_{\ast}(u) \, \vert \, \mathrm{Q}\big)  \, \mathrm{d}u.
\end{equation}
\begin{proof} Let $s,s_{0} \geqslant 0$ and choose $T \geqslant \max\{s,s_{0}\}$. Taking expectations under $\mathds{P}^{\boldsymbol{\gamma}_{\ast}}$ in \hyperref[loc.m.c]{(\ref*{loc.m.c})}, and noting that the stochastic integral process in \hyperref[loc.m.b]{(\ref*{loc.m.b})} is a martingale, leads to
\begin{equation} \label{cor.a.a}
\mathds{E}_{\mathds{P}^{\boldsymbol{\gamma}_{\ast}}}\big[\Lambda\big(s,\overline{X}(s)\big)\big] - \mathds{E}_{\mathds{P}^{\boldsymbol{\gamma}_{\ast}}}\big[ \Lambda\big(s_{0},\overline{X}(s_{0})\big)\big]
= - \tfrac{1}{2} \int_{s_{0}}^{s}  \mathds{E}_{\mathds{P}^{\boldsymbol{\gamma}_{\ast}}}\Big[ \big\vert \nabla \Lambda\big( u,\overline{X}(u)\big) \big\vert^{2}\Big]  \, \mathrm{d}u.
\end{equation}
Recalling the relative entropy \hyperref[rel.ent.sec]{(\ref*{rel.ent.sec})} and the relative Fisher information \hyperref[fish.inf]{(\ref*{fish.inf})}, we obtain \hyperref[cor.a.b]{(\ref*{cor.a.b})}.
\end{proof}
\end{corollary}

\begin{corollary} \label{cor.b} For Lebesgue-almost every $s \geqslant 0$, the rate of relative entropy dissipation equals
\begin{equation}
\frac{\mathrm{d}}{\mathrm{d}s} \, H\big( \overline{P}_{\ast}(s) \, \vert \, \mathrm{Q} \big)  
= - \tfrac{1}{2} \, I\big( \overline{P}_{\ast}(s) \, \vert \, \mathrm{Q}\big).
\end{equation}
\end{corollary}


\section{From backwards dynamics ``back'' to forward dynamics} \label{yet.st.cont.}


Starting with the forward Langevin--Smoluchowski dynamics \hyperref[ls.dp.sde]{(\ref*{ls.dp.sde})}, we have seen in \hyperref[sec.ascp.f]{Section \ref*{sec.ascp.f}} that the combined effects of time-reversal, and of stochastic control of the drift under an entropic-type criterion, lead to the backwards dynamics
\begin{equation} \label{blsd.rev}
\mathrm{d} \overline{X}(s) =  - \nabla \Psi\big(\overline{X}(s)\big)  \, \mathrm{d}s + \mathrm{d}\overline{W}^{\boldsymbol{\gamma}_{\ast}}(s), \qquad 0 \leqslant s \leqslant T,
\end{equation}
which are again of the Langevin--Smoluchowski type, but now viewed on the filtered probability space $(\Omega,\overline{\mathds{G}},\mathds{P}^{\boldsymbol{\gamma}_{\ast}})$. We will see now that this universal property of Langevin--Smoluchowski measure is consistent in the following sense: starting with the backwards Langevin--Smoluchowski dynamics of \hyperref[blsd.rev]{(\ref*{blsd.rev})}, after another reversal of time, the solution of a related stochastic control problem leads to the original forward Langevin--Smoluchowski dynamics \hyperref[ls.dp.sde]{(\ref*{ls.dp.sde})} we started with. This consistency property should come as no surprise, but its formal proof requires the results of \hyperref[sec.ltm]{Section \ref*{sec.ltm}}, which perhaps appeared artificial at first sight. 

\smallskip

Let us recall from \hyperref[rev.time.again]{Subsection \ref*{rev.time.again}} that reversing time in \hyperref[blsd.rev]{(\ref*{blsd.rev})} leads to the forward dynamics 
\begin{equation} \label{forward.dyn.as}
\mathrm{d} X(t) = \Big( \nabla \Lambda\big(T-t,X(t)\big) - \nabla \Psi\big(X(t)\big) \Big) \, \mathrm{d}t + \mathrm{d}W^{\boldsymbol{\gamma}_{\ast}}(t), \qquad 0 \leqslant t \leqslant T
\end{equation}
on the filtered probability space $(\Omega,\mathds{F},\mathds{P}^{\boldsymbol{\gamma}_{\ast}})$. By analogy with \hyperref[sec.ascp.f]{Section \ref*{sec.ascp.f}}, we define an equivalent probability measure $\Pi^{\boldsymbol{\beta}} \sim \mathds{P}^{\boldsymbol{\gamma}_{\ast}}$ as follows.

\smallskip

For any measurable process $[0,T] \times \Omega \ni (t,\omega) \mapsto \boldsymbol{\beta}(t,\omega) \in \mathds{R}^{n}$, adapted to the forward filtration $\mathds{F}$ of \hyperref[can.forw.filt]{(\ref*{can.forw.filt})} and satisfying the condition
\begin{equation} \label{adap.pro.sec.three}
\mathds{P}^{\boldsymbol{\gamma}_{\ast}}\bigg[ \int_{0}^{T} \vert \boldsymbol{\beta}(t)\vert^{2} \, \mathrm{d}t < \infty \bigg] = 1,
\end{equation}  
we consider the exponential $(\mathds{F},\mathds{P}^{\boldsymbol{\gamma}_{\ast}})$-local martingale
\begin{equation} \label{exp.mart.two.three}
Z^{\boldsymbol{\beta}}(t) \vcentcolon = \exp \bigg( \int_{0}^{t} \big\langle \boldsymbol{\beta}(\theta) \, , \, \mathrm{d}W^{\boldsymbol{\gamma}_{\ast}}(\theta) \big\rangle - \tfrac{1}{2} \int_{0}^{t} \vert \boldsymbol{\beta}(\theta)\vert^{2} \, \mathrm{d}\theta \bigg) 
\, , \qquad 0 \leqslant t \leqslant T.
\end{equation}
We denote by $\mathcal{B}$ the collection of all processes $\boldsymbol{\beta}$ as above, for which $Z^{\boldsymbol{\beta}}$ is a true $(\mathds{F},\mathds{P}^{\boldsymbol{\gamma}_{\ast}})$-martingale. 

\smallskip

Now, for every $\boldsymbol{\beta} \in \mathcal{B}$, we introduce an equivalent probability measure $\Pi^{\boldsymbol{\beta}} \sim \mathds{P}^{\boldsymbol{\gamma}_{\ast}}$ on path space, via
\begin{equation} \label{girs.dens.sec.three}
\frac{\mathrm{d}\Pi^{\boldsymbol{\beta}}}{\mathrm{d}\mathds{P}^{\boldsymbol{\gamma}_{\ast}}} \bigg\vert_{\mathcal{F}(t)} 
= Z^{\boldsymbol{\beta}}(t)
\, , \qquad 0 \leqslant t \leqslant T.
\end{equation}
Then we deduce from the Girsanov theorem that, under the probability measure $\Pi^{\boldsymbol{\beta}}$, the process
\begin{equation} \label{ti.re.ga.bm.sec.three}
W^{\boldsymbol{\beta}}(t) \vcentcolon = W^{\boldsymbol{\gamma}_{\ast}}(t) - \int_{0}^{t} \boldsymbol{\beta}(\theta) \, \mathrm{d}\theta, \qquad 0 \leqslant t \leqslant T
\end{equation}
is $\mathds{R}^{n}$-valued $\mathds{F}$-Brownian motion, and the dynamics \hyperref[forward.dyn.as]{(\ref*{forward.dyn.as})} become
\begin{equation} \label{ls.dp.sde.sec.b.three}
\mathrm{d} X(t) = \Big( \nabla \Lambda\big(T-t,X(t)\big) + \boldsymbol{\beta}(t) - \nabla \Psi\big(X(t)\big) \Big) \, \mathrm{d}t + \mathrm{d}W^{\boldsymbol{\beta}}(t), \qquad 0 \leqslant t \leqslant T.
\end{equation}

\smallskip

We couple these dynamics with the stochastic differential \hyperref[loc.m.c]{(\ref*{loc.m.c})} and deduce that the process
\begin{equation} \label{for.sub.mar}
N^{\boldsymbol{\beta}}(t) \vcentcolon = \Lambda\big(T-t,X(t)\big) + \tfrac{1}{2} \int_{0}^{t} \vert \boldsymbol{\beta}(\theta)\vert^{2} \, \mathrm{d}\theta, \qquad 0 \leqslant t \leqslant T
\end{equation}
is a local $\Pi^{\boldsymbol{\beta}}$-submartingale with decomposition 
\begin{equation} \label{fif.de}
\mathrm{d}N^{\boldsymbol{\beta}}(t) 
= \tfrac{1}{2} \big\vert \nabla \Lambda\big(T-t,X(t)\big) + \boldsymbol{\beta}(t) \big\vert^{2}  \, \mathrm{d}t + \Big \langle \nabla \Lambda\big(T-t,X(t)\big) \, , \, \mathrm{d}\overline{W}^{\boldsymbol{\beta}}(t) \Big\rangle.
\end{equation}
In fact, introducing the sequence
\begin{equation} \label{stop.time.four}
\tau_{n} \vcentcolon = \inf \bigg\{ t \geqslant 0 \colon \, \int_{0}^{t} \Big( \big\vert \nabla \Lambda\big(T-\theta,X(\theta)\big) \big\vert^{2}  + \vert \boldsymbol{\beta}(\theta)\vert^{2} \Big) \, \mathrm{d}\theta \, \geqslant n \,  \bigg\} \wedge T
\, , \qquad n \in \mathds{N}_{0}
\end{equation}
of $\mathds{F}$-stopping times with $\tau_{n} \uparrow T$, we see that the stopped process $N^{\boldsymbol{\beta}}(\, \cdot \, \wedge \tau_{n})$ is an $\mathds{F}$-submartingale under $\Pi^{\boldsymbol{\beta}}$, for every $n \in \mathds{N}_{0}$. In particular, we observe
\begingroup
\addtolength{\jot}{0.7em}
\begin{equation} \label{tr.st.subm..four}
\begin{aligned} 
H\big( P(2T) \, \vert \, \mathrm{Q} \big) &= H\big( \overline{P}_{\ast}(T) \, \vert \, \mathrm{Q} \big) = \mathds{E}_{\mathds{P}^{\boldsymbol{\gamma}_{\ast}}}\big[ \Lambda\big(T,\overline{X}(T)\big) \big] = \mathds{E}_{\Pi^{\boldsymbol{\beta}}}\big[ \Lambda\big(T,X(0)\big) \big] \\
&\leqslant \mathds{E}_{\Pi^{\boldsymbol{\beta}}}\bigg[ \Lambda\big(T-\tau_{n},X(\tau_{n})\big) + \tfrac{1}{2} \int_{0}^{\tau_{n}} \vert \boldsymbol{\beta}(\theta) \vert^{2} \, \mathrm{d}\theta \bigg],
\end{aligned}
\end{equation}
\endgroup
since $\Pi^{\boldsymbol{\beta}} = \mathds{P}^{\boldsymbol{\gamma}_{\ast}}$ on $\mathcal{F}(0) = \boldsymbol{\sigma}(X(0)) = \boldsymbol{\sigma}(\overline{X}(T))$. Now \hyperref[tr.st.subm..four]{(\ref*{tr.st.subm..four})} holds for every $n \in \mathds{N}_{0}$, thus
\begin{equation} \label{tr.st.subm.lim.four}
H\big( P(2T) \, \vert \, \mathrm{Q} \big) 
\leqslant \liminf_{n \rightarrow \infty} \ \mathds{E}_{\Pi^{\boldsymbol{\beta}}}\bigg[ \Lambda\big(T-\tau_{n},X(\tau_{n})\big) + \tfrac{1}{2} \int_{0}^{\tau_{n}} \vert \boldsymbol{\beta}(\theta) \vert^{2} \, \mathrm{d}\theta \bigg].
\end{equation}

\smallskip

But the minimum in \hyperref[sl.st.pde.m.sec.three]{(\ref*{sl.st.pde.m.sec.three})} is attained by $b_{\ast} = - \nabla \Lambda(s,x)$; likewise, the drift term in \hyperref[fif.de]{(\ref*{fif.de})} vanishes, if we select the process $\boldsymbol{\beta}_{\ast} \in \mathcal{B}$ via
\begin{equation} \label{score.pro.four}
\boldsymbol{\beta}_{\ast}(t,\omega) \vcentcolon = - \nabla \Lambda\big(T-t,\omega(t)\big), \qquad \textnormal{ thus } \qquad \boldsymbol{\beta}_{\ast}(t)  = - \nabla \Lambda\big(T-t,X(t)\big), \qquad 0 \leqslant t \leqslant T.
\end{equation}
With this choice, the forward dynamics of \hyperref[ls.dp.sde.sec.b.three]{(\ref*{ls.dp.sde.sec.b.three})} take the form
\begin{equation} \label{bwd.ls.sec.thir.four}
\mathrm{d} X(t) = - \nabla \Psi\big(X(t)\big) \, \mathrm{d}t + \mathrm{d}W^{\boldsymbol{\beta}_{\ast}}(t);
\end{equation}
that is, \textit{precisely the forward Langevin--Smoluchowski dynamics} \hyperref[ls.dp.sde]{(\ref*{ls.dp.sde})} we started with, but now with the ``initial condition'' $X(0) = \overline{X}(T)$ and independent driving $\mathds{F}$-Brownian motion $W^{\boldsymbol{\beta}_{\ast}}$, under $\Pi^{\boldsymbol{\beta}_{\ast}}$. Since $\Pi^{\boldsymbol{\beta}_{\ast}} = \mathds{P}^{\boldsymbol{\gamma}_{\ast}}$ holds on the $\boldsymbol{\sigma}$-algebra $\mathcal{F}(0) = \boldsymbol{\sigma}(X(0)) = \boldsymbol{\sigma}(\overline{X}(T))$, the initial distribution of $X(0)$ under $\Pi^{\boldsymbol{\beta}_{\ast}}$ is equal to $P(2T)$. Furthermore, with $\boldsymbol{\beta} = \boldsymbol{\beta}_{\ast}$, the process of \hyperref[for.sub.mar]{(\ref*{for.sub.mar})}, \hyperref[fif.de]{(\ref*{fif.de})} becomes a $\Pi^{\boldsymbol{\beta}_{\ast}}$-local martingale, namely
\begin{equation} \label{loc.mart.four}
N^{\boldsymbol{\beta}_{\ast}}(t) 
= \Lambda\big(T,X(0)\big) + \int_{0}^{t} \Big \langle \nabla \Lambda\big(T-\theta,X(\theta)\big) \, , \, \mathrm{d}\overline{W}^{\boldsymbol{\beta}_{\ast}}(\theta) \Big\rangle, \qquad 0 \leqslant t \leqslant T;
\end{equation}
and we have equality in \hyperref[tr.st.subm..four]{(\ref*{tr.st.subm..four})}, thus also
\begin{equation} \label{tr.st.subm.lim.sec.four}
H\big( P(2T) \, \vert \, \mathrm{Q} \big) 
= \lim_{n \rightarrow \infty} \ \mathds{E}_{\Pi^{\boldsymbol{\beta}_{\ast}}}\bigg[ \Lambda\big(T-\tau_{n},X(\tau_{n})\big) + \tfrac{1}{2} \int_{0}^{\tau_{n}} \vert \boldsymbol{\beta}_{\ast}(\theta) \vert^{2} \, \mathrm{d}\theta \bigg].
\end{equation}

\smallskip

We conclude that the infimum over $\boldsymbol{\beta} \in \mathcal{B}$ of the right-hand side in \hyperref[tr.st.subm.lim.four]{(\ref*{tr.st.subm.lim.four})} is attained by the process $\boldsymbol{\beta}_{\ast}$ of \hyperref[score.pro.four]{(\ref*{score.pro.four})}, which gives rise to the Langevin--Smoluchowski dynamics \hyperref[bwd.ls.sec.thir.four]{(\ref*{bwd.ls.sec.thir.four})} for the process $(X(t))_{0 \leqslant t \leqslant T}$, under $\Pi^{\boldsymbol{\beta}_{\ast}}$. We formalize this result as follows.

\begin{theorem} \label{m.theo.one.four} Consider the stochastic control problem of minimizing over the class $\mathcal{B}$ of measurable, adapted processes $\boldsymbol{\beta}$ satisfying \textnormal{\hyperref[adap.pro.sec.three]{(\ref*{adap.pro.sec.three})}} and inducing an exponential martingale $Z^{\boldsymbol{\beta}}$ in \textnormal{\hyperref[exp.mart.two.three]{(\ref*{exp.mart.two.three})}}, with the notation of \textnormal{\hyperref[stop.time.four]{(\ref*{stop.time.four})}} and with the forward dynamics of \textnormal{\hyperref[ls.dp.sde.sec.b.three]{(\ref*{ls.dp.sde.sec.b.three})}}, the expected cost 
\begin{equation} \label{expected.cost.four}
\mathcal{J}(\boldsymbol{\beta}) \vcentcolon = \liminf_{n \rightarrow \infty} \ \mathds{E}_{\Pi^{\boldsymbol{\beta}}}\bigg[ \Lambda\big(T-\tau_{n},X(\tau_{n})\big) + \tfrac{1}{2} \int_{0}^{\tau_{n}} \vert \boldsymbol{\beta}(\theta) \vert^{2} \, \mathrm{d}\theta \bigg].
\end{equation}

\smallskip

Under the assumptions of \textnormal{\hyperref[the.setting]{Section \ref*{the.setting}}}, the infimum $\inf_{\boldsymbol{\beta} \in \mathcal{B}} \mathcal{J}(\boldsymbol{\beta})$ is equal to the relative entropy $H( P(2T) \, \vert \, \mathrm{Q})$ and is attained by the ``score process'' $\boldsymbol{\beta}_{\ast}$ of \textnormal{\hyperref[score.pro.four]{(\ref*{score.pro.four})}}. This choice leads to the forward Langevin--Smoluchowski dynamics \textnormal{\hyperref[bwd.ls.sec.thir.four]{(\ref*{bwd.ls.sec.thir.four})}}, and with $\boldsymbol{\beta} = \boldsymbol{\beta}_{\ast}$ the limit in \textnormal{\hyperref[expected.cost.four]{(\ref*{expected.cost.four})}} exists as in \textnormal{\hyperref[tr.st.subm.lim.sec.four]{(\ref*{tr.st.subm.lim.sec.four})}}. 
\begin{proof} We have to show that the minimizing process $\boldsymbol{\beta}_{\ast}$ belongs to the collection $\mathcal{B}$ of admissible processes. By its definition in \hyperref[score.pro.four]{(\ref*{score.pro.four})}, the process $\boldsymbol{\beta}_{\ast}$ is measurable, and adapted to the forward filtration $\mathds{F}$ of \hyperref[can.forw.filt]{(\ref*{can.forw.filt})}. Thanks to \hyperref[true.exp.a]{(\ref*{true.exp.a})} in \hyperref[rel.ent.prop.]{Proposition \ref*{rel.ent.prop.}}, we have 
\begin{equation} \label{finite.kst.four}
\mathds{E}_{\mathds{P}^{\boldsymbol{\gamma}_{\ast}}}\bigg[ \int_{0}^{T} \vert  \boldsymbol{\beta}_{\ast}(t) \vert^{2} \, \mathrm{d}t\bigg]
= \mathds{E}_{\mathds{P}^{\boldsymbol{\gamma}_{\ast}}}\bigg[ \int_{0}^{T} \big\vert \nabla \Lambda\big( T-t,X(t)\big) \big\vert^{2} \, \mathrm{d}t\bigg] < \infty,
\end{equation}
which implies a fortiori that the condition in \hyperref[adap.pro.sec.three]{(\ref*{adap.pro.sec.three})} is satisfied for $\boldsymbol{\beta} = \boldsymbol{\beta}_{\ast}$.

\smallskip

It remains to check that the process $Z^{\boldsymbol{\beta}_{\ast}}$ defined in the manner of \hyperref[exp.mart.two.three]{(\ref*{exp.mart.two.three})}, is a true martingale. From \hyperref[rel.ent.prop.]{Proposition \ref*{rel.ent.prop.}} we have the stochastic differential
\begin{equation}
\mathrm{d}\Lambda\big(T-t,X(t)\big) 
= \tfrac{1}{2} \big\vert \nabla \Lambda\big( T-t,X(t)\big) \big\vert^{2}  \, \mathrm{d}t + \Big \langle \nabla \Lambda\big(T-t,X(t)\big) \, , \, \mathrm{d}W^{\boldsymbol{\gamma}_{\ast}}(t) \Big\rangle,
\end{equation}
and therefore
\begingroup
\addtolength{\jot}{0.7em}
\begin{align}
&\int_{0}^{t} \big\langle \boldsymbol{\beta}_{\ast}(\theta) \, , \, \mathrm{d}W^{\boldsymbol{\gamma}_{\ast}}(\theta) \big\rangle - \tfrac{1}{2} \int_{0}^{t} \vert \boldsymbol{\beta}_{\ast}(\theta)\vert^{2} \, \mathrm{d}\theta \\
& \qquad = - \int_{0}^{t} \Big\langle \nabla \Lambda\big(T-\theta,X(\theta)\big) \, , \, \mathrm{d}W^{\boldsymbol{\gamma}_{\ast}}(\theta) \Big\rangle - \tfrac{1}{2} \int_{0}^{t} \big\vert \nabla \Lambda\big(T-\theta,X(\theta)\big) \big\vert^{2} \, \mathrm{d}\theta \\
& \qquad = \Lambda\big(T,X(0)\big) - \Lambda\big(T-t,X(t)\big) = \log \bigg( \frac{\lambda\big(T,X(0)\big)}{\lambda\big(T-t,X(t)\big)} \bigg),
\end{align}
\endgroup
which expresses the exponential process of \hyperref[exp.mart.two.three]{(\ref*{exp.mart.two.three})} with $\boldsymbol{\beta} = \boldsymbol{\beta}_{\ast}$ as
\begin{equation}
Z^{\boldsymbol{\beta}_{\ast}}(t) = \frac{\lambda\big(T,X(0)\big)}{\lambda\big(T-t,X(t)\big)} \, , \qquad 0 \leqslant t \leqslant T.
\end{equation}

\smallskip

The process $Z^{\boldsymbol{\beta}_{\ast}}$ is a positive local martingale, thus a supermartingale. To see that it is a true $(\mathds{F},\mathds{P}^{\boldsymbol{\gamma}_{\ast}})$-martingale, it suffices to argue that it has constant expectation. But $Z^{\boldsymbol{\beta}_{\ast}}(0) \equiv 1$, so we have to show $\mathds{E}_{\mathds{P}^{\boldsymbol{\gamma}_{\ast}}}[Z^{\boldsymbol{\beta}_{\ast}}(T)] = 1$. We denote again by $P(s,y;t,\xi)$ the transition kernel of the Langevin--Smoluchowski dynamics, note that 
\begin{equation}
\mathds{P}^{\boldsymbol{\gamma}_{\ast}}\big[\overline{X}(s) \in \mathrm{d}y, \overline{X}(t) \in \mathrm{d}\xi\big] = \overline{p}_{\ast}(s,y) \, P(s,y;t,\xi) \, \mathrm{d}y \, \mathrm{d}\xi
\end{equation}
for $0 \leqslant s < t \leqslant T$ and $(y,\xi) \in \mathds{R}^{n} \times \mathds{R}^{n}$, and recall the invariance property \hyperref[invar.gam]{(\ref*{invar.gam})} of $\mathrm{Q}$, to deduce $\mathds{E}_{\mathds{P}^{\boldsymbol{\gamma}_{\ast}}}[Z^{\boldsymbol{\beta}_{\ast}}(T)] = 1$ in the manner of \hyperref[invar.gam.a]{(\ref*{invar.gam.a})} -- \hyperref[invar.gam.d]{(\ref*{invar.gam.d})}.
This implies that $Z^{\boldsymbol{\beta}_{\ast}}$ is a true martingale and completes the proof of \hyperref[m.theo.one.four]{Theorem \ref*{m.theo.one.four}}.
\end{proof}
\end{theorem}


\subsection{Discussion: entropic interpretation of the expected cost when \texorpdfstring{$\mathrm{Q}(\mathds{R}^{n}) < \infty$}{Q is a probability measure}} \label{int.two.four}


By analogy with \hyperref[int.one]{Subsection \ref*{int.one}}, we interpret now the expected cost on the right-hand side of \hyperref[expected.cost.four]{(\ref*{expected.cost.four})} in terms of relative entropies. From \hyperref[exp.mart.two.three]{(\ref*{exp.mart.two.three})} -- \hyperref[ti.re.ga.bm.sec.three]{(\ref*{ti.re.ga.bm.sec.three})}, we deduce that
\begin{equation} \label{log.girs.dens.rel.ent.t.four}  
\log \bigg(\frac{\mathrm{d}\Pi^{\boldsymbol{\beta}}}{\mathrm{d}\mathds{P}^{\boldsymbol{\gamma}_{\ast}}} \bigg\vert_{\mathcal{F}(\tau_{n})} \bigg)
=  \int_{0}^{\tau_{n}} \big\langle \boldsymbol{\beta}(\theta) \, , \, \mathrm{d}W^{\boldsymbol{\beta}}(\theta) \big\rangle + \tfrac{1}{2} \int_{0}^{\tau_{n}} \vert \boldsymbol{\beta}(\theta)\vert^{2} \, \mathrm{d}\theta  
\end{equation}
holds for every $\boldsymbol{\beta} \in \mathcal{B}$ and $n \in \mathds{N}_{0}$. Therefore, as the $\Pi^{\boldsymbol{\beta}}$-expectation of the stochastic integral in \hyperref[log.girs.dens.rel.ent.t.four]{(\ref*{log.girs.dens.rel.ent.t.four})} vanishes, the expected quadratic cost, or ``energy'', term in \hyperref[expected.cost.four]{(\ref*{expected.cost.four})} is equal to the relative entropy
\begin{equation}
\mathds{E}_{\Pi^{\boldsymbol{\beta}}}\bigg[ \tfrac{1}{2} \int_{0}^{\tau_{n}} \vert \boldsymbol{\beta}(\theta) \vert^{2} \, \mathrm{d}\theta \bigg]
= \mathds{E}_{\Pi^{\boldsymbol{\beta}}}\bigg[ \log \bigg(  \frac{\mathrm{d}\Pi^{\boldsymbol{\beta}}}{\mathrm{d}\mathds{P}^{\boldsymbol{\gamma}_{\ast}}} \bigg\vert_{\mathcal{F}(\tau_{n})} \bigg)\bigg].
\end{equation}

\smallskip

In order to interpret the first term in \hyperref[expected.cost.four]{(\ref*{expected.cost.four})}, let us assume that $\mathrm{Q}$, and thus also the induced measure $\mathds{Q}$ on path space, are probability measures. Then we have
\begingroup
\addtolength{\jot}{0.7em}
\begin{align}
\mathds{E}_{\Pi^{\boldsymbol{\beta}}}\big[ \Lambda\big(T-\tau_{n},X(\tau_{n})\big) \big]
&= \mathds{E}_{\Pi^{\boldsymbol{\beta}}}\bigg[ \log \bigg( \frac{\mathrm{d}\mathds{P}^{\boldsymbol{\gamma}_{\ast}}}{\mathrm{d}\mathds{Q}} \bigg\vert_{\boldsymbol{\sigma}(X(\tau_{n}))} \bigg)\bigg] \label{narrent.at.four} \\
&= \mathds{E}_{\Pi^{\boldsymbol{\beta}}}\bigg[ \log \bigg( \frac{\mathrm{d}\Pi^{\boldsymbol{\beta}}}{\mathrm{d}\mathds{Q}} \bigg\vert_{\boldsymbol{\sigma}(X(\tau_{n}))} \bigg)\bigg] 
- \mathds{E}_{\Pi^{\boldsymbol{\beta}}}\bigg[ \log \bigg( \frac{\mathrm{d}\Pi^{\boldsymbol{\beta}}}{\mathrm{d}\mathds{P}^{\boldsymbol{\gamma}_{\ast}}} \bigg\vert_{\boldsymbol{\sigma}(X(\tau_{n}))} \bigg)\bigg]. \label{narrent.bt.four}
\end{align}
\endgroup
We conclude that, in this case, the expected cost of \hyperref[expected.cost.four]{(\ref*{expected.cost.four})} is equal to the sum $H_{n}^{\boldsymbol{\beta}} + D_{n}^{\boldsymbol{\beta}}$ of two non-negative quantities:
\begin{equation}
H_{n}^{\boldsymbol{\beta}} \vcentcolon = \mathds{E}_{\Pi^{\boldsymbol{\beta}}}\bigg[ \log \bigg( \frac{\mathrm{d}\Pi^{\boldsymbol{\beta}}}{\mathrm{d}\mathds{Q}} \bigg\vert_{\boldsymbol{\sigma}(X(\tau_{n}))} \bigg)\bigg] ,
\end{equation}
the relative entropy of the probability measure $\Pi^{\boldsymbol{\beta}}$ with respect to the probability measure $\mathds{Q}$ when both are restricted to the $\boldsymbol{\sigma}$-algebra generated by the random variable $X(\tau_{n})$; and
\begin{equation} \label{ent.difft.four}
D_{n}^{\boldsymbol{\beta}} \vcentcolon 
= \mathds{E}_{\Pi^{\boldsymbol{\beta}}}\bigg[ \log \bigg( \frac{\mathrm{d}\Pi^{\boldsymbol{\beta}}}{\mathrm{d}\mathds{P}^{\boldsymbol{\gamma}_{\ast}}} \bigg\vert_{\mathcal{F}(\tau_{n})} \bigg)\bigg]
- \mathds{E}_{\Pi^{\boldsymbol{\beta}}}\bigg[ \log \bigg( \frac{\mathrm{d}\Pi^{\boldsymbol{\beta}}}{\mathrm{d}\mathds{P}^{\boldsymbol{\gamma}_{\ast}}} \bigg\vert_{\boldsymbol{\sigma}(X(\tau_{n}))} \bigg)\bigg],
\end{equation}
the difference between the relative entropies of the probability measure $\Pi^{\boldsymbol{\beta}}$ with respect to the probability measure $\mathds{P}^{\boldsymbol{\gamma}_{\ast}}$, when restricted to the $\boldsymbol{\sigma}$-algebra generated by the collection of random variables $(X(\theta \wedge \tau_{n}))_{0 \leqslant \theta \leqslant T}$ and by the random variable $X(\tau_{n})$, respectively. 


\section{The case of finite invariant measure, and an iterative procedure} \label{conclusion.seven}


Let us suppose now that the diffusion process $(X(t))_{t \geqslant 0}$ as in \hyperref[ls.dp.sde]{(\ref*{ls.dp.sde})} is well-defined, along with the curve $P(t) = \mathrm{Law}(X(t))$, $t \geqslant 0$ of its time-marginals; and that the invariant measure $\mathrm{Q}$ of \hyperref[subsec.2.1]{Subsection \ref*{subsec.2.1}} is finite, i.e., \hyperref[prob.gib]{(\ref*{prob.gib})} holds, and is thus normalized to a probability measure.

\smallskip

Then, neither the coercivity condition \hyperref[docc]{(\ref*{docc})}, nor the finite second moment condition \hyperref[int.req.a]{(\ref*{int.req.a})}, are needed for the results of \hyperref[sec.ascp.f]{Sections \ref*{sec.ascp.f}} -- \hyperref[yet.st.cont.]{\ref*{yet.st.cont.}}. The reason is that the relative entropy $H(P(t) \, \vert \, \mathrm{Q})$ is now well-defined and non-negative, as both $P(t)$ and $\mathrm{Q}$ are probability measures. Since the function $t \mapsto H(P(t) \, \vert \, \mathrm{Q})$ is decreasing and the initial relative entropy $H(P(0) \, \vert \, \mathrm{Q})$ is finite on account of \hyperref[int.req.b]{(\ref*{int.req.b})}, it follows that this function takes values in $[0,\infty)$. It can also be shown in this case that 
\begin{equation} \label{ent.dec.lim}
\lim_{t \rightarrow \infty} \downarrow H\big(P(t) \, \vert \, \mathrm{Q}\big) = 0,
\end{equation}
i.e., the relative entropy decreases down to zero; see \cite[Proposition 1.9]{FJ16} for a quite general version of this result. This, in turn, implies that the time-marginals $(P(t))_{t \geqslant 0}$ converge to $\mathrm{Q}$ in total variation as $t \rightarrow \infty$, on account of the \textit{Pinsker--Csisz\'{a}r inequality} \begin{equation} 
2 \, \Vert P(t) - \mathrm{Q}\Vert_{\mathrm{TV}}^{2} \leqslant H\big( P(t) \, \vert \, \mathrm{Q}\big).
\end{equation}
The entropic decrease to zero is actually exponentially fast, whenever the Hessian of the potential $\Psi$ dominates a positive multiple of the identity matrix; see, e.g., \cite{BE85}, \cite[Section 5]{MV00}, \cite[Proposition 1']{OV00}, \cite[Formal Corollary 9.3]{Vil03}, or \cite[Remark 3.23]{KST20b}. As another consequence of \hyperref[ent.dec.lim]{(\ref*{ent.dec.lim})}, the initial relative entropy $H(P(0) \, \vert \, \mathrm{Q})$ can be expressed as
\begin{equation} \label{i.r.e.i}
H\big( P(0) \, \vert \, \mathrm{Q} \big) = \tfrac{1}{2} \, \mathds{E}_{\mathds{P}}\bigg[ \int_{0}^{\infty} \big \vert \nabla L\big(t,X(t)\big) \big\vert^{2} \, \mathrm{d}t\bigg]. 
\end{equation}
We prove \hyperref[ent.dec.lim]{(\ref*{ent.dec.lim})} and \hyperref[i.r.e.i]{(\ref*{i.r.e.i})} in \hyperref[app.a.a]{Appendix \ref*{app.a.a}}.

\smallskip
 
In this context, i.e., with \hyperref[prob.gib]{(\ref*{prob.gib})} replacing \hyperref[docc]{(\ref*{docc})} and \hyperref[int.req.a]{(\ref*{int.req.a})}, and always under the standing assumption \hyperref[int.req.b]{(\ref*{int.req.b})}, Theorem 4.1 in \cite{KST20a} continues to hold, as do the results in \hyperref[sec.ascp.f]{Sections \ref*{sec.ascp.f}} -- \hyperref[yet.st.cont.]{\ref*{yet.st.cont.}}. By combining time-reversal with stochastic control of the drift, these results lead to an alternating sequence of forward and backward Langevin--Smoluchowski dynamics, with time-marginals starting at $P(0)$ and converging along $(P(kT))_{k \in \mathds{N}_{0}}$ in total variation to the invariant probability measure $\mathrm{Q}$. Along the way, the values of the corresponding stochastic control problems decrease along $(H(P(kT) \, \vert \, \mathrm{Q}))_{k \in \mathds{N}}$ to zero.


\setkomafont{section}{\large}
\setkomafont{subsection}{\normalsize}


\begin{appendices}

\section{The decrease of the relative entropy without convexity assumption} \label{app.a.a}

We present a probabilistic proof of \hyperref[ent.dec.lim]{(\ref*{ent.dec.lim})} and \hyperref[i.r.e.i]{(\ref*{i.r.e.i})}, which complements the proof of the more general Proposition 1.9 in \cite{FJ16}. We stress that no convexity assumptions are imposed on the potential $\Psi$.
 
\begin{proof}[Proof of \texorpdfstring{\textnormal{\hyperref[ent.dec.lim]{(\ref*{ent.dec.lim})}}}{} and \texorpdfstring{\textnormal{\hyperref[i.r.e.i]{(\ref*{i.r.e.i})}}}{}:] Since $\mathrm{Q}$ is assumed to be a probability measure in \hyperref[conclusion.seven]{Section \ref*{conclusion.seven}}, the relative entropy $H(P(t) \, \vert \, \mathrm{Q})$ is non-negative for every $t \geqslant 0$. Thus, \cite[Corollary 4.3]{KST20a} gives the inequality
\begin{equation} \label{7.2.a}
H\big( P(0) \, \vert \, \mathrm{Q} \big) = H\big( P(T) \, \vert \, \mathrm{Q} \big) + \tfrac{1}{2} \, \mathds{E}_{\mathds{P}}\bigg[ \int_{0}^{T} \big \vert \nabla L\big(t,X(t)\big) \big\vert^{2} \, \mathrm{d}t\bigg] 
\geqslant \tfrac{1}{2} \, \mathds{E}_{\mathds{P}}\bigg[ \int_{0}^{T} \big \vert \nabla L\big(t,X(t)\big) \big\vert^{2} \, \mathrm{d}t\bigg]
\end{equation}
for every $T \in (0,\infty)$. Letting $T \uparrow \infty$ in \hyperref[7.2.a]{(\ref*{7.2.a})}, we deduce from the monotone convergence theorem that
\begin{equation} \label{7.2.b}
\tfrac{1}{2} \, \mathds{E}_{\mathds{P}}\bigg[ \int_{0}^{\infty} \big \vert \nabla L\big(t,X(t)\big) \big\vert^{2} \, \mathrm{d}t\bigg] \leqslant H\big( P(0) \, \vert \, \mathrm{Q} \big).
\end{equation}

\smallskip

By analogy with \hyperref[t.p.s]{Subsection \ref*{t.p.s}}, we denote by $\mathds{Q}$ the Langevin--Smoluchowski measure associated with the potential $\Psi$, but now with distribution 
\begin{equation} \label{7.2.c}
\mathds{Q}\big[ X(0) \in A \big] = \mathrm{Q}(A) = \int_{A} q(x) \, \mathrm{d}x, \qquad A \in \mathscr{B}(\mathds{R}^{n})
\end{equation}
for the random variable $X(0)$. Since $\mathrm{Q}$ is a probability measure, the Langevin--Smoluchowski measure $\mathds{Q}$ is a well-defined probability measure on path space $\Omega = C([0,\infty);\mathds{R}^{n})$. 

\smallskip

Let us recall now the likelihood ratio of \hyperref[l.r.f]{(\ref*{l.r.f})}. We denote the corresponding likelihood ratio process by $\vartheta(t) \vcentcolon = \ell(t,X(t))$, $t \geqslant 0$. The following remarkable insight is due to Pavon \cite{Pav89}, and was rediscovered by Fontbona and Jourdain \cite{FJ16}: For any given $T \in (0,\infty)$, the time-reversed likelihood ratio process
\begin{equation} \label{7.2.d}
\overline{\vartheta}(s) \vcentcolon = \ell\big(T-s,\overline{X}(s)\big) =  \frac{p\big(T-s,X(T-s)\big)}{q\big(\overline{X}(s)\big)}  \, , \qquad 0 \leqslant s \leqslant T
\end{equation}
is a $\mathds{Q}$-martingale with respect to the backwards filtration $\overline{\mathds{G}} = (\overline{\mathcal{G}}(s))_{0 \leqslant s \leqslant T}$ of \hyperref[can.pro]{(\ref*{can.pro})}. For a simple proof of this result in the setting of this paper we refer to \cite[Appendix E]{KST20b}. 

\smallskip

Let us pick arbitrary times $0 \leqslant t_{1} < t_{2} < \infty$. For any given $T \in (t_{2},\infty)$, the martingale property of the process \hyperref[7.2.d]{(\ref*{7.2.d})} amounts, with $s_{1} = T-t_{1}$ and $s_{2} = T-t_{2}$, to
\begin{equation} \label{7.2.e}
\mathds{E}_{\mathds{Q}}\Big[ \overline{\vartheta}(s_{1}) \, \big\vert \, \overline{\mathcal{G}}(s_{2})\Big] = \overline{\vartheta}(s_{2})
\qquad \Longleftrightarrow \qquad
\mathds{E}_{\mathds{Q}}\Big[ \vartheta(t_{1}) \, \big\vert \, \boldsymbol{\sigma}\big( X(\theta) \colon t_{2} \leqslant \theta \leqslant T \big) \Big] = \vartheta(t_{2}).
\end{equation}
Because $T \in (t_{2},\infty)$ is arbitrary in \hyperref[7.2.e]{(\ref*{7.2.e})}, this gives
\begin{equation} \label{7.2.f}
\mathds{E}_{\mathds{Q}}\big[ \vartheta(t_{1}) \, \vert \, \mathcal{H}(t_{2}) \big] = \vartheta(t_{2}) 
\, , \qquad 
\mathcal{H}(t) \vcentcolon = \boldsymbol{\sigma}\big( X(\theta) \colon t \leqslant \theta < \infty \big).
\end{equation}
In other words, the likelihood ratio process $(\vartheta(t))_{t \geqslant 0}$ is a backwards $\mathds{Q}$-martingale of the backwards filtration $\mathds{H} = (\mathcal{H}(t))_{t \geqslant 0}$. We denote by $\mathcal{H}(\infty) \vcentcolon = \bigcap_{t \geqslant 0} \mathcal{H}(t)$ the tail $\boldsymbol{\sigma}$-algebra of the Langevin--Smoluchowski diffusion $(X(t))_{t \geqslant 0}$. The ergodicity of this process under the probability measure $\mathds{Q}$ implies that the tail $\boldsymbol{\sigma}$-algebra $\mathcal{H}(\infty)$ is $\mathds{Q}$-trivial, i.e., $\mathcal{H}(\infty) = \{ \varnothing, \Omega \}$ modulo $\mathds{Q}$; see \hyperref[app.a]{Appendix \ref*{app.a}} for a proof of this claim.

\smallskip

We recall now the martingale version of Theorem 9.4.7 (backwards submartingale convergence) in \cite{Chu01}. This says that $(\vartheta(t))_{t \geqslant 0}$ is a $\mathds{Q}$-uniformly integrable family, that the limit
\begin{equation} \label{7.2.g}
\vartheta(\infty) \vcentcolon = \lim_{t \rightarrow \infty} \vartheta(t)
\end{equation}
exists $\mathds{Q}$-a.e., that the convergence in \hyperref[7.2.g]{(\ref*{7.2.g})} holds also in $L^{1}(\mathds{Q})$, and that for every $t \geqslant 0$ we have
\begin{equation} \label{7.2.h}
\mathds{E}_{\mathds{Q}}\big[ \vartheta(t) \, \vert \, \mathcal{H}(\infty) \big] = \vartheta(\infty) \, , \qquad \mathds{Q}\textnormal{-a.e.}
\end{equation}
But since the tail $\boldsymbol{\sigma}$-algebra $\mathcal{H}(\infty)$ is $\mathds{Q}$-trivial, the random variable $\vartheta(\infty)$ is $\mathds{Q}$-a.e.\ constant, and \hyperref[7.2.h]{(\ref*{7.2.h})} identifies this constant as $\vartheta(\infty) \equiv 1$.

\smallskip

In terms of the function $f(x) \vcentcolon = x \log x$ for $x > 0$ (and with $f(0) \vcentcolon = 0$), we can express the relative entropy $H(P(t) \, \vert \, \mathrm{Q})$ as
\begin{equation} \label{7.2.i}
H\big( P(t) \, \vert \, \mathrm{Q} \big) 
= \mathds{E}_{\mathds{P}}\big[ \log \vartheta(t) \big]
= \mathds{E}_{\mathds{Q}}\big[ f\big(\vartheta(t)\big) \big] \geqslant 0 \, , \qquad t \geqslant 0.
\end{equation}
The convexity of $f$, in conjunction with \hyperref[7.2.f]{(\ref*{7.2.f})}, shows that the process $(f(\vartheta(t)))_{t \geqslant 0}$ is a backwards $\mathds{Q}$-submartingale of the backwards filtration $\mathds{H}$, with decreasing expectation as in \hyperref[7.2.i]{(\ref*{7.2.i})}. By appealing to the backwards submartingale convergence theorem \cite[Theorem 9.4.7]{Chu01} once again, we deduce that $(f(\vartheta(t)))_{t \geqslant 0}$ is a $\mathds{Q}$-uniformly integrable family, which converges, a.e.\ and in $L^{1}$ under $\mathds{Q}$, to
\begin{equation} \label{7.2.j}
\lim_{t \rightarrow \infty} f\big(\vartheta(t)\big) = f\big(\vartheta(\infty)\big) = f(1) = 0.
\end{equation}
In particular,
\begin{equation} \label{7.2.k}
\lim_{t \rightarrow \infty} \downarrow H\big(P(t) \, \vert \, \mathrm{Q}\big) 
= \lim_{t \rightarrow \infty} \mathds{E}_{\mathds{Q}}\big[ f\big(\vartheta(t)\big) \big]
= \mathds{E}_{\mathds{Q}}\big[  \lim_{t \rightarrow \infty}  f\big(\vartheta(t)\big) \big]
=0,
\end{equation}
proving \hyperref[ent.dec.lim]{(\ref*{ent.dec.lim})}. From \hyperref[7.2.k]{(\ref*{7.2.k})} and \hyperref[7.2.a]{(\ref*{7.2.a})} it follows now that \hyperref[7.2.b]{(\ref*{7.2.b})} holds as equality, proving \hyperref[i.r.e.i]{(\ref*{i.r.e.i})}.
\end{proof}

\section{The triviality of the tail \texorpdfstring{$\boldsymbol{\sigma}$}{sigma}-algebra \texorpdfstring{$\mathcal{H}(\infty)$}{}} \label{app.a}

We recall the filtered probability space $(\Omega,\mathcal{F}(\infty),\mathds{F},\mathds{Q})$. Here, $\Omega = C([0,\infty);\mathds{R}^{n})$ is the path space of continuous functions, $\mathcal{F}(\infty) = \boldsymbol{\sigma}(\bigcup_{t \geqslant 0} \mathcal{F}(t))$, the canonical filtration $\mathds{F} = (\mathcal{F}(t))_{t \geqslant 0}$ is as in \hyperref[can.forw.filt]{(\ref*{can.forw.filt})}, and the Langevin--Smoluchowski measure $\mathds{Q}$ is represented by
\begin{equation} \label{a.1}
\mathds{Q}(B) = \int_{\mathds{R}^{n}} \mathds{P}^{x}(B) \, \mathrm{d}\mathrm{Q}(x) \, , \qquad B \in \mathscr{B}(\Omega),
\end{equation}
where $\mathds{P}^{x}$ denotes the Langevin--Smoluchowski measure with initial distribution $\delta_{x}$, for every $x \in \mathds{R}^{n}$, and $\mathscr{B}(\Omega)$ is the Borel $\boldsymbol{\sigma}$-field\footnote{There are several equivalent ways to define the Borel $\boldsymbol{\sigma}$-field on $\Omega$. Two possible constructions appear in Problems 2.4.1 and 2.4.2 in \cite{KS88}.} on $\Omega$.

\smallskip

For every $s \geqslant 0$, we define a measurable map $\theta_{s} \colon \Omega \rightarrow \Omega$, called \textit{shift transformation}, by requiring that $\theta_{s}(\omega)(t) = \omega(s+t)$ holds for all $\omega \in \Omega$ and $t \geqslant 0$. A Borel set $B \in \mathscr{B}(\Omega)$ is called \textit{shift-invariant} if $\theta_{s}^{-1}(B) = B$ holds for any $s \geqslant 0$. Since the Gibbs probability measure $\mathrm{Q}$ is the unique invariant measure for the Langevin--Smoluchowski diffusion $(X(t))_{t \geqslant 0}$, \cite[Theorem 3.8]{Bel06} implies that the probability measure $\mathds{Q}$ of \hyperref[a.1]{(\ref*{a.1})} is \textit{ergodic}, meaning that $\mathds{Q}(B) \in \{0,1\}$ holds for every shift-invariant set $B$. As a consequence of the ergodicity of $\mathds{Q}$, the \textit{Birkhoff Ergodic Theorem} \cite[Theorem 3.4]{Bel06} implies that, for every $A \in \mathscr{B}(\mathds{R}^{n})$, the limit
\begin{equation} \label{a.2}
\lim_{t \rightarrow \infty} \frac{1}{t} \int_{0}^{t} \mathds{1}_{A}\big(X(s)\big) \, \mathrm{d}s = \mathrm{Q}(A)
\end{equation}
exists $\mathds{Q}$-a.e.

\begin{lemma} The tail $\boldsymbol{\sigma}$-algebra $\mathcal{H}(\infty)$ is $\mathds{Q}$-trivial, i.e., $\mathcal{H}(\infty) = \{ \varnothing, \Omega \}$ modulo $\mathds{Q}$.
\begin{proof} We follow a reasoning similar to that in \cite[Remark 1.10]{FJ16}. According to \cite[Theorem 1.3.9]{Kun90}, it suffices to show that the Langevin--Smoluchowski diffusion $(X(t))_{t \geqslant 0}$ is \textit{recurrent in the sense of Harris}, i.e.,
\begin{equation} \label{a.3}
\mathds{P}^{x} \bigg[ \int_{0}^{\infty} \mathds{1}_{A}\big(X(s)\big) \, \mathrm{d}s = \infty \bigg] = 1 
\end{equation}
is satisfied for every $x \in \mathds{R}^{n}$ and all $A \in \mathscr{B}(\mathds{R}^{n})$ with $\mathrm{Q}(A) > 0$.

\smallskip

For the proof of \hyperref[a.3]{(\ref*{a.3})}, we fix $x \in \mathds{R}^{n}$ and $A \in \mathscr{B}(\mathds{R}^{n})$ with $\mathrm{Q}(A) > 0$. By its definition, the event 
\begin{equation} \label{a.4}
B \vcentcolon = \bigg\{ \int_{0}^{\infty} \mathds{1}_{A}\big(X(s)\big) \, \mathrm{d}s = \infty \bigg\} 
\end{equation}
is shift-invariant. Thus, by the ergodicity of $\mathds{Q}$, the probability $\mathds{Q}(B)$ is equal to either zero or one. Outside the set $B$, we have
\begin{equation} \label{a.5}
\lim_{t \rightarrow \infty} \frac{1}{t} \int_{0}^{t} \mathds{1}_{A}\big(X(s)\big) \, \mathrm{d}s = 0.
\end{equation}
But since $\mathrm{Q}(A) > 0$, the $\mathds{Q}$-a.e.\ limit \hyperref[a.2]{(\ref*{a.2})} implies that $\mathds{Q}(B^{\mathrm{c}}) = 0$ and hence $\mathds{Q}(B) = 1$. From the definition \hyperref[a.1]{(\ref*{a.1})} of the probability measure $\mathds{Q}$ it follows that $\mathds{P}^{x}(B) = 1$ for $\mathrm{Q}$-a.e.\ $x \in \mathds{R}^{n}$. Since $\mathrm{Q}$ is equivalent to Lebesgue measure, we also have that $\mathds{P}^{x}(B) = 1$ for Lebesgue-a.e.\ $x \in \mathds{R}^{n}$.

\smallskip

Furthermore, the shift-invariance of $B$ and the Markov property of the Langevin--Smoluchowski diffusion give
\begin{equation} \label{a.6}
\mathds{P}^{x}(B) = \mathds{P}^{x} \big( \theta_{t}^{-1}(B) \big) = \mathds{E}_{\mathds{P}^{x}}\big[ \mathds{P}^{X(t)}(B) \big] = T_{t}\big( \mathds{P}^{x}(B)\big) = \int_{\mathds{R}^{n}} P(0,x;t,\mathrm{d}y) \, \mathds{P}^{y}(B)
\end{equation}
for every $t \geqslant 0$. Here, $P(0,x;t,y)$ denotes the transition kernel of the Langevin--Smoluchowski dynamics, so that $\mathds{P}^{x}[X(t) \in \mathrm{d}y] = P(0,x;t,y) \, \mathrm{d}y$; and $T_{t}$ denotes the operator 
\begin{equation}
T_{t}f(x) \vcentcolon = \int_{\mathds{R}^{n}} P(0,x;t,\mathrm{d}y) \, f(y) \, , \qquad (t,x) \in [0,\infty) \times \mathds{R}^{n}    
\end{equation}
acting on bounded measurable functions $f \colon \mathds{R}^{n} \rightarrow \mathds{R}$. Since $(T_{t})_{t \geqslant 0}$ is a strong Feller semigroup under the assumptions of this paper, the function $\mathds{R}^{n} \ni x \mapsto T_{t}f(x)$ is continuous. Now \hyperref[a.6]{(\ref*{a.6})} implies the continuity of the function $\mathds{R}^{n} \ni x \mapsto \mathds{P}^{x}(B)$. On the other hand, we have already seen that the function $\mathds{R}^{n} \ni x \mapsto \mathds{P}^{x}(B) \in [0,1]$ is Lebesgue-a.e.\ equal to one. But such a function is constant everywhere, i.e., $\mathds{P}^{x}(B) = 1$ for every $x \in \mathds{R}^{n}$, proving \hyperref[a.3]{(\ref*{a.3})}.
\end{proof}
\end{lemma}
\end{appendices}




\bibliographystyle{alpha}
{\footnotesize
\bibliography{references}}


\end{document}